\documentclass[12pt]{article}

\usepackage{geometry}
\usepackage{hyperref} 

\usepackage{amsmath}
\usepackage{amssymb}
\usepackage{amscd}
\usepackage{textcomp}
\usepackage{dsfont}
\usepackage{mathrsfs}

\long\def\forget#1{}

\newcommand{\Verkuerzung}[2]{#1}

\usepackage{color}
\newcounter{commentcounter}
\newcommand{\comment}[1]{\stepcounter{commentcounter}{\color{red}\textbf{Comment \arabic{commentcounter}.} #1}
\immediate\write16{}
\immediate\write16{Warning: There was still a comment . . . }
\immediate\write16{}}

\newcounter{urscommentcounter}

\def\?{\ 
{\bf\color{red}???}\ 
\immediate\write16{}
\immediate\write16{Warning: There was still a question mark . . . }
\immediate\write16{}}

\usepackage{amsthm}
\theoremstyle{plain}
\newtheorem{theorem}{Theorem}[section]

\newtheorem{lemma}[theorem]{Lemma}

\newtheorem{corollary}[theorem]{Corollary}
\newtheorem{proposition}[theorem]{Proposition}
    
\newtheorem{question}[theorem]{Question}

\theoremstyle{definition}
\newtheorem{definition}[theorem]{Definition}

\newtheorem{definition-theorem}[theorem]{Definition-Theorem}
\newtheorem{definition-remark}[theorem]{Definition-Remark}
\newtheorem{notation-remark}[theorem]{Notation-Remark}

\newtheorem{remark}[theorem]{Remark}

\theoremstyle{remark}

\usepackage{xy}
\xyoption{all}

\newdir^{ (}{{}*!/-3pt/\dir^{(}}    
\newdir^{  }{{}*!/-3pt/\dir^{}}    
\newdir_{ (}{{}*!/-3pt/\dir_{(}}    
\newdir_{  }{{}*!/-3pt/\dir_{}}

\newcounter{zahl}

\def\theenumi{(\alph{enumi})}

\def\p@enumii{\theenumi}

\newcommand{\DS}{\displaystyle}

\newcommand{\SC}{\scriptstyle}

\DeclareMathOperator{\Aut}{Aut}

\DeclareMathOperator{\End}{End}

\DeclareMathOperator{\GL}{GL}

\DeclareMathOperator{\Koh}{H}

\DeclareMathOperator{\Hom}{Hom}

\DeclareMathOperator{\Isom}{Isom}

\DeclareMathOperator{\Pic}{Pic}

\DeclareMathOperator{\SL}{SL}

\DeclareMathOperator{\Spec}{Spec}
\DeclareMathOperator{\Spf}{Spf}

\DeclareMathOperator{\Tr}{Tr}

\DeclareMathOperator{\charakt}{char}

\DeclareMathOperator{\diag}{diag}

\DeclareMathOperator{\equi}{equi}
\newcommand{\et}{{\acute{e}t\/}}

\DeclareMathOperator{\id}{\,id}

\renewcommand{\phi}{\varphi}
\renewcommand{\epsilon}{\varepsilon}

\usepackage{amsfonts}
\newcommand{\BOne} {{\mathchoice{\hbox{\rm1\kern-2.7pt l\kern.9pt}}
                              {\hbox{\rm1\kern-2.7pt l\kern.9pt}}
                              {\hbox{\scriptsize\rm1\kern-2.3pt l\kern.4pt}}
                              {\hbox{\scriptsize\rm1\kern-2.4pt l\kern.5pt}}}}

\newcommand{\BA}{{\mathbb{A}}}

\newcommand{\BD}{{\mathbb{D}}}

\newcommand{\BF}{{\mathbb{F}}}
\newcommand{\BG}{{\mathbb{G}}}

\newcommand{\BP}{{\mathbb{P}}}
\newcommand{\BQ}{{\mathbb{Q}}}

\newcommand{\BZ}{{\mathbb{Z}}}

\newcommand{\CA}{{\cal{A}}}

\newcommand{\CC}{{\cal{C}}}

\newcommand{\CF}{{\cal{F}}}
\newcommand{\CG}{{\cal{G}}}
\newcommand{\CH}{{\cal{H}}}

\newcommand{\CL}{{\cal{L}}}

\newcommand{\CN}{{\cal{N}}}
\newcommand{\CO}{{\cal{O}}}
\newcommand{\CP}{{\cal{P}}}
\newcommand{\CQ}{{\cal{Q}}}

\newcommand{\CS}{{\cal{S}}}
\newcommand{\CT}{{\cal{T}}}
\newcommand{\CU}{{\cal{U}}}
\newcommand{\CV}{{\cal{V}}}
\newcommand{\CW}{{\cal{W}}}
\newcommand{\CX}{{\cal{X}}}
\newcommand{\CY}{{\cal{Y}}}
\newcommand{\CZ}{{\cal{Z}}}

\newcommand{\FF}{{\mathfrak{F}}}
\newcommand{\FG}{{\mathfrak{G}}}

\newcommand{\FX}{{\mathfrak{X}}}

\newcommand{\FZ}{{\mathfrak{Z}}}

\newcommand{\Fm}{{\mathfrak{m}}}

\newcommand{\scrH}{{\mathscr{H}}}

\let\setminus\smallsetminus

\newcommand{\ul}[1]{{\underline{#1}}}
\newcommand{\ol}[1]{{\overline{#1}}}
\newcommand{\wh}[1]{{\widehat{#1}}}
\newcommand{\wt}[1]{{\widetilde{#1}}}

\usepackage{ifthen}

\newcommand{\invlim}[1][]{\ifthenelse{\equal{#1}{}}
{\DS \lim_{\longleftarrow}}
{\DS \lim_{\underset{#1}{\longleftarrow}}}
}

\newcommand{\dirlim}[1][]{\ifthenelse{\equal{#1}{}}
{\DS \lim_{\longrightarrow}}
{\DS \lim_{\underset{#1}{\longrightarrow}}}
}

\newcommand{\dotBD}{\vbox{\hbox{\kern2pt\bf.}\vskip-4.5pt\hbox{$\BD$}}}

\DeclareMathOperator{\Nilp}{\CN \!{\it ilp}}

\DeclareMathOperator{\Sets}{\CS \!{\it ets}}

\def\longto{\longrightarrow}
\def\into{\hookrightarrow}

\def\isoto{\stackrel{}{\mbox{\hspace{1mm}\raisebox{+1.4mm}{$\SC\sim$}\hspace{-3.5mm}$\longrightarrow$}}}

\newbox\mybox
\def\arrover#1{\mathrel{
       \setbox\mybox=\hbox spread 1.4em{\hfil$\scriptstyle#1$\hfil}
       \vbox{\offinterlineskip\copy\mybox
             \hbox to\wd\mybox{\rightarrowfill}}}}

\newcommand{\genericG}{P}

\newcommand{\tauGlob}{\tau}

\newcommand{\charsect}{s}

\newcommand{\scrE}{\mathscr{E}}

\begin{document}

\author{Esmail Arasteh Rad\forget{\footnote{?}} }

\date{\today}

\title{A Note On\\ Two Fiber Bundles and
The Manifestations Of\\ ``Shtuka''\\}

\maketitle

\begin{abstract}

In this note we intend to look at the moduli stacks for global $\FG$-shtukas from a new perspective. We discuss a unifying interpretation of several moduli spaces (stacks) including \forget{(mixed characteristic) }moduli of global $\FG$-shtukas and (a variant of the) moduli of Higgs bundles. We view these spaces (stacks) as different fibers of a family over a scheme (stack) locally of finite type. We discuss (a relative version of) the local model theory for this family. We also consider the Hecke stacks over the moduli stack of $\FG$-shtukas and discuss the corresponding (motivic) Hecke operations. 

\end{abstract}

\tableofcontents

\section*{Introduction}\label{SectIntroduction}

Let $C$ be a smooth projective geometrically irreducible curve over a perfect field $k$. Let $\FG$ be a flat affine group scheme of finite type over $C$. The stack $\scrH^1(C,\FG)$ of (principal) $\FG$-bundles over $C$, has been considered by various authors, especially due to the key roles that it plays, both in the geometric Langlands program and also in the arithmetic Langlands program over function fields. This stack is an Artin stack, locally of finite type; e.g. see \cite{Beh} when $\FG$ is constant split reductive, and \cite{AH_Global} for the general case. \\
 
There are interesting although inadequate studies aimed to describe the symmetries of $\scrH^1(C,\FG)$. Note that when $\FG=\GL_r$, this moduli stack coincides the moduli of vector bundles of rank $r$ over $C$. In $\charakt k=0$, and in the presence of stability condition, in \cite{BGM} the authors describe the group of automorphisms of the moduli space of stable vector bundles. They observe that they all come from the obvious ones. This means that they either arise from $\Aut(C)$, or by twisting by a line bundle, or possibly by sending a vector bundle to its dual vector bundle. Note in addition that the stack of morphisms between Artin stacks has been considered, and studied, by several authors, including Aoki, Olson, Hall and Rydh, see \cite{Ao}, \cite{Ols2} and \cite{HR}.\\

On the other hand, for some purposes, it is useful to consider the moduli of $\FG$-bundles over a relative curve $C$ over $S$. In \cite{Wang}, Wang proves that \forget{for affine group scheme $\FG$ over a base field $k$, }the stack $\scrH^1(C,\FG)$ remains Artin (over $S$) when we replace the curve $C/k$ with a projective scheme $X$ over a base $S$. Considering the relative version was in fact motivated and proposed for certain applications in the geometric Langlands program. Beyond this, as another remarkable example, to prove the purity of the cohomology of this moduli stack over $\BF_p$, in \cite{HeiSch}, Heinloth and Schmitt, consider the relative situation over a curve $C$ over a Dedekind domain $R$ over $\BZ$, and implement techniques from nearby-vanishing cycles to lift the Atiyah-Bott theory from characteristic zero (in the generic fiber) to the special fiber of this moduli stack.\\

In this note we introduce a family $\Sigma$ over the stack $\scrE$ of endomorphisms of $\scrH^1(C,\FG)$, and we study the geometry of the family and its fibers. Note that after imposing relevant boundedness conditions (and endowing with extra structures), different fibers of this family realize some interesting moduli stacks (spaces), such as (global Schubert varieties inside) Beilinson-Drinfeld affine Grassmannian, (a variant of the) moduli of Higgs bundles, and moduli of (bounded) $\FG$-shtukas; see subsection \ref{SubsectTwoFiberBundles} and subsection \ref{Subsec_SpecificFibers}. \\

Let us briefly explain the content of each section. In Section \ref{SectPrelim}, after we fix some notation in \ref{SubsectNotation}, we recall the basic definitions of formal algebraic stacks (and rigid analytic stacks) in Subsection \ref{SubsectFormalAndAnalyticStacks}. In Section \ref{SectGeneralConstruction}, we explain the general construction of the family $\Sigma=\Sigma(\CH\rightrightarrows \CY)\to \scrE$, corresponding to a two fiber bundle $\CH\rightrightarrows\CY$. Here $\scrE$ is the stack of endomorphisms of $\CY$. We prove a general local model theorem for a fiber over an endomorphism of $\CY$ which annihilate the tangent bundle; see Theorem \ref{ThmLocModTFB}. Then we discuss some applications of this theorem; see Corollaries \ref{CorICSheafNablaH} and \ref{CorSingularities}.


\noindent
In paragraph \ref{SubsectionHeinlothSchmitt} we recall Heinloth-Schmitt stability condition. We will observe that at least for split reductive case, the stack of endomorphisms of $\scrH^1(C,\FG)$, after imposing stability condition, is an Artin stack. In paragraph \ref{SubsectBDGrassmannian} we recall that the Beilinson-Drinfeld affine Grassmannian also arises in this context, and we further discuss boundedness conditions and functoriality. We further discuss how a boundedness condition gives rise to a Hecke cycle; see Proposition \ref{PropHeckOverBunG}.\\ 
\noindent
In paragraph \ref{SubsectHiggsBundles} we consider certain symmetries of the moduli stack $\scrH^1(C,\FG)$ which are encapsulated in itself, in the sense of twisting by torsors. Consequently we see that certain fibers can be regarded as a variant of the moduli stacks of Higgs bundles. We further discuss the (formally) properness of certain restrictions of this family, as well as some lifting properties of the corresponding stacks. 

\noindent
In paragraph \ref{SectModuliG-Sht} we discuss the corresponding picture for the moduli of $\FG$-shtukas. We address similar lifting problems. See Theorem \ref{Thm_FormalDM} and Proposition \ref{PropKohSpec&generic}. \\
In Section \ref{Sect The Hecke stack over the moduli of G-Shtukas}, as another example for a two fiber bundle over an algebraic stack, we introduce a Hecke stack over moduli of $\FG$-shtukas. We in particular observe that loop group invariant cycles inside a global affine Grassmannian $GR_m(C,\FG)$ induce certain homomorphisms between motives of the moduli stacks of $\FG$-shtukas; see theorem \ref{ThmCyclesOnNablaH}. The later gives certain bivariant classes. 

\section*{Acknowledgment}
I would like to thank Urs Hartl and Chia-Fu Yu for their encouragement and continued support, as well as many helpful conversations.\\
I have benefited and inspired by deep philosophical ideas of Arash Rastegar and wish to express my deep appreciation to him. I would like to thank Eaman Eftekhary for interesting conversations and very useful comments related to certain structures arising from two fiber bundles. I also thank Giuseppe Ancona and Somayeh Habibi for some helpful explanations related to the theory of motives.\\
\noindent
This manuscript has been extracted out of the notes of my lecture in S\'eminaire de Th\'eorie des Nombres at Caen and for that I would like to thank the organizers of the seminar. I would especially like to thank Tuan Ng\^o Dac and warmly thank him for very interesting and insightful conversations, as well as his kind and efficient support.
 
\section{Preliminaries}\label{SectPrelim}

\subsection{Notation and Conventions}\label{SubsectNotation}

\begin{tabbing}

$\genericG_\nu:=\FG\times_C\Spec Q_\nu,$\; \=\kill

$S$\> a locally noetherian scheme,\\[1mm]

$C$\> a smooth projective relative curve over $S$ with geometrically\\  ~~~~~~~~~~~~~~~~~~~~~~~~~~~~~irreducible fibers,\\[1mm]


$\FG$\> a smooth affine group scheme over $C$,\\[1mm]

$Sch/S$\> category of schemes over $S$,\\[1mm]

$n$\> a positive integer,\\[1mm]

$\ul s$\> an $n$-tuple of sections $s_i:T\to C$,\\[1mm]

$\Gamma_\ul s$\> the union of the graphs $\Gamma_{s_i}\subseteq C\times_S T$,\\[1mm]




\end{tabbing}

\bigskip
Assume that we have two morphisms $f,g\colon X\to Y$ of schemes or stacks. We denote by $\equi(f,g\colon X\rightrightarrows Y)$ the pull back of the diagonal under the morphism $(f,g)\colon X\to Y\times_\BZ Y$, that is 
$$
\equi(f,g\colon X\rightrightarrows Y)\,:=\,X\times_{(f,g),Y\times Y,\Delta}Y
$$ 
where $\Delta=\Delta_{Y/\BZ}\colon Y\to Y\times_\BZ Y$ is the diagonal morphism.\\

\bigskip

By an IC-sheaf $IC(\CX)$ on a stack $\CX$, we will mean the intermediate extension of the constant perverse sheaf $\ol \BQ_\ell$ on an open dense substack $\CX^\circ$ of $\CX$ such that the corresponding reduced stack $\CX_{red}^\circ$ is smooth. The IC-sheaf is normalized so that it is pure of weight zero.

\bigskip

Let $\wh S$ be a formal scheme. We denote by $\Nilp_\wh S$ the category of schemes over $\wh S$ on which an ideal of definition of  $\wh S$ is locally nilpotent, equipped with the \'etale topology. 

\bigskip

Let $H$ be a sheaf of groups (for the \'etale topology) on a scheme $X$. By a (\emph{right}) \emph{$H$-torsor} (also called an \emph{$H$-bundle}) on $X$ we mean a sheaf $\CG$ for the \'etale topology on $X$ together with a (right) action of the sheaf $H$ such that $\CG$ is isomorphic to $H$ on a \'etale covering of $X$. Here $H$ is viewed as an $H$-torsor by right multiplication. We denote by $\scrH^1(C,\FG)$ the category fibered in groupoids over $Sch/S$ with $\scrH^1(C,\FG)(T)$ the groupoid of $\FG$-bundles over $C_T:=C\times_S T$.

\bigskip

When $S=\Spec \BF_q$ we follow our notation in \cite{AH_LM} and \cite{AH_Global}, in particular by $\nabla_n\scrH^1(C,\FG)$ we denote the moduli stack whose $T$-points consists of $\FG$-shtukas over $T$ with $n$ characteristic sections, e.g. see \cite[Definition 3.3]{AH_Global} or \cite[2.0.10]{AH_LM}. 

\bigskip

For a perfect field $k$ and $X$ in $\textbf{Sch}_k$, let $Ch_i(X)$ and $Ch^i(X)$ denote Fulton's $i$-th Chow groups and let $Ch_\ast(X):=\oplus_i Ch_i(X)$ (resp. $Ch^\ast(X):=\oplus_i Ch^i(X)$).  \\

\noindent
Finally, to denote the motivic categories over a perfect field $k$, such as 

$$
\textbf{DM}_{gm}(k,\BZ),~ \textbf{DM}_{gm}^{eff}(k,\BZ),~ \textbf{DM}_{-}^{eff}(k),~\text{etc.} 
$$  
and the functors $M(-):\textbf{Sch}_k\rightarrow \textbf{DM}_{gm}^{eff}(k,\BZ)$ and $M^c(-):\textbf{Sch}_k\rightarrow \textbf{DM}_{gm}^{eff}(k,\BZ)$ we use the same notation that was introduced  in \cite{VSF}. 
When $char~k>0$ we assume coefficients in $\BQ$. 
\\

For the definition of the geometric motives with compact support in positive characteristic we also refer to \cite[Appendix B]{H-K}.\\

\subsection{Formal and analytic stacks}\label{SubsectFormalAndAnalyticStacks}

Recall that a \emph{formal space} $\wh X$ over a formal scheme $\wh S$ is a sheaf of sets on the site $\Nilp_{\wh S}$. In addition it is called a \emph{formal algebraic space} if the diagonal morphism $\wh X \to \wh X \times_{\wh S} \wh X$ is representable by a quasi-compact
morphism of formal schemes and there is a formal scheme $\wh X'$ over $\wh S$ and a morphism of formal $\wh S$-spaces $\wh X'\to \wh X$ which is representable by an \'etale surjective morphism
of formal schemes.\\

Let $\CX$ be a stack over a scheme $S$. Let $S_0$ be a locally closed subscheme of $S$. Let $\wh S$ denote the formal completion of $S$ along $S_0$. Restricting the fibered functor $\CX$ to the category $\Nilp_{\wh S}$ gives a category $\wh \CX$ fibred in groupoids over $\Nilp_{\wh S}$ which inherits the following properties from $\CX$ 

\begin{enumerate}
\item[i)] for every $\CV$ in $\Nilp_{\wh S}$ and $x, y$ in $\wh{\CX}(\CV)$ the presheaf
\begin{eqnarray*}
\Isom: Sch\slash \CV &\to &  \Sets \\
\CU\to \CV &\longmapsto & \Hom_{\wh{\CX}(\CU)}(x_\CU, y_\CU),
\end{eqnarray*}

is a sheaf on $Sch\slash \CV$.

\item[ii)] for every covering $\CV_i \to \CV$ in $\Nilp_{\wh S}$ all descent data for this covering are effective.

\end{enumerate}

Furthermore if $\CX$ is an algebraic stack in the sense of Artin (resp. Deligne-Mumford (DM)) we have

\begin{enumerate}

\item

the diagonal 1-morphism $\wh \CX \to \wh \CX \times_{\wh S} \wh \CX$  over $\wh S$
is representable by formal algebraic $\wh S$-spaces (resp. schemes),\forget{ i.e. the fiber over any $\CU$-valued point in $\Nilp_{\wh S}$ is representable by a formal algebraic $\wh S$-space (resp. $\wh S$-scheme),} separated, and quasi-compact.

\forget{
the diagonal 1-morphism $\wh \CX \to \wh \CX \times_{\wh S} \wh \CX$  over $\wh S$
is representable (i.e. the fiber over any $\CU$-valued point in $\Nilp_{\wh S}$ is representable by a formal algebraic $\wh S$-space), separated, and quasi-compact,
\comment{representable by formal algebraic $\wh S$-spaces (resp. schemes), i.e. the fiber over any $\CU$-valued point in $\Nilp_{\wh S}$ is representable by a formal algebraic $\wh S$-space (resp. $\wh S$-scheme), separated, and quasi-compact}
}
\item
there exists a formal algebraic $\wh S$-space $\wh X$ and a presentation 
$$
P : \wh X \to \wh \CX
$$
 of formal $\wh S$-stacks which is representable by a smooth (resp. \'etale) and surjective morphism of formal algebraic $\wh S$-spaces. 

\end{enumerate}

\noindent 
Let us abstractify the above observation and phrase it in the following way

\begin{definition}
A category $\wh \CX$ fibered in groupoids over $\Nilp_{\wh S}$ is called a \emph{formal stack} if it has the properties i) and ii) indicated above. Also we say $\wh \CX$ is \emph{formal algebraic (Artin) stack} (resp. formal stack of DM-type) if in addition it is subject to a) and b) above.   
\end{definition}

Let $K$ be a complete discrete valuation field. Let $\CO_K$ and $k$ denote the corresponding ring of integers and residue field.

\begin{remark}\label{RemAnalyticStack}
Let $\CA n_K$ denote the category of
$K$-analytic spaces, equipped with the \'etale topology. A $K$-analytic stack $\CX$ is a stack in groupoids over the site $\CA n_K$. It is called Artin (resp. Deligne Mumford) if similar conditions to the above conditions (a) and (b) hold in this category. One can extend the analytification functor  from the category of schemes to obtain the functor $(-)^{an}$ from the category of algebraic stacks locally of finite type over $K$ to the category of $K$-analytic stacks. We
refer to \cite[Section~6]{PY} for details. Similarly, one can produce the special fiber functor $(-)_s$  (resp. the generic fiber
functor $(-)_\eta$) from
the category of formal stacks locally finitely presented over $\Spf \CO_K$ to the category of algebraic stacks locally of finite type over $k$ (resp. to the category of $K$-analytic stacks). 

\end{remark}

\section{Two fiber bundles and the dynamics of the base}\label{SectGeneralConstruction}

In this section we introduce a family $\Sigma$, corresponding to a two fiber bundle $\CH$, over certain stacks of endomorphisms. We prove some general statements related to local and global geometry of certain fibers of such families. 

\subsection{General construction of the family $\Sigma\to \scrE$}\label{SubsectTwoFiberBundles}

\begin{definition}\label{Def_HR}
\begin{enumerate}

 Let $\CY$ (resp. $\CY'$) be an algebraic stack in the Artin's sense (resp. formal algebraic stack), locally of finite type over a scheme $S$. Let $X$ be a projective flat scheme over $S$.\\
\item Let $\CH$ be a stack over $X$ via a morphism $char:\CH\to X$, and furthermore assume that $\CH$ is fibred over $\CY$ and $\CY'$ via the following maps

$$
pr^{\leftarrow}: \CH \to \CY ~\text{and}~ pr^{\rightarrow}: \CH \to \CY'.
$$

We call the tuple $\ul\CH:=(\CH,~char, pr^\leftarrow,pr^\rightarrow)$, consisting of the above data, a $HR$-tuple. We say that $\ul\CH$ is a $CHR$-tuple of degree $m$ if $pr^\rightarrow\times char$ is proper and of finite type, and $pr^\leftarrow$ is finite type and flat of relative dimension $m$. 

\item
Let $\ul\Hom:=\Hom(\CY,\CY')$ denote the corresponding $\Hom$-stack. It is contravariant 2-functor from the category of affine noetherian schemes over $S$ to the 2-category of groupoids given by assigning the groupoid of 1-morphisms $\Hom_T(\CY\times_S T,\CY'\times_S T)$ to a test scheme $T$. When $\CY'=\CY$, then we set $\scrE:=\ul\Hom(\CY,\CY)$.

\item
Define $\Sigma(\ul\CH)$ via the following pull-back diagram
$$
\CD
\Sigma(\ul\CH)@>>>\CH\\
@V{\phi}VV @VV(pr^{\leftarrow}, pr^{\rightarrow})V\\
\ul\Hom \times \CY'@>>>\CY\times \CY'
\endCD
$$ 
of stacks.
\noindent
The bottom arrow is given by
$$
(f:\CY \to \CY',y) \mapsto (y, f(y)).
$$

We view $\Sigma(\ul\CH)$ as a family over $\ul\Hom$. When it is clear from the context we use the shorthand $\Sigma$ to denote  $\Sigma(\ul\CH)$. 

\end{enumerate}
\end{definition}

\begin{remark}

Let $\ul\CH$ be a $CHR$-tuple of degree m. Then one may define the following Hecke-type operation, from the category of perverse sheaves on $\CY$ to the derived category of sheaves on $\CY'\times_S X$ given by
the formula

$$
(pr^\rightarrow\times char)_{\ast} \circ pr^{\leftarrow\ast} (-).
$$
\end{remark}

\begin{notation-remark}
One can mimic the above construction also in the category of formal stacks over a formal scheme $\wh S$.  We then use the notation $\wh \CY$, $\wh \CY'$, $\wh\CH$, $\ul{\wh\CH}$, $\wh\scrE$, $\wh\Sigma:=\wh\Sigma(\wh{\ul\CH})$ and etc. to denote the corresponding formal stacks.
\end{notation-remark}

\bigskip
\begin{remark}\label{Remark_HilbStack}
Assume that $\CY$ is a quasi-projective scheme over $S$. Then there is an obvious morphism $\scrE\to Hilb_{\dim \CY}(\CY\times_S\CY)$, defined by sending $f$ to its graph $\Gamma_f\subseteq \CY\times_S \CY$. This identifies $\scrE$ with an open subscheme of $Hilb_{\dim \CY}(\CY\times_S\CY)$. In particular each connected component of $\scrE$ is of finite type, and these components form a countable set. The automorphism group scheme $Aut(\CY)$ is open in $\scrE$ by \cite[p. 267]{Gro} (see also \cite[Lemma I.1.10.1]{Kol}). If $\CY$ is a projective variety, then $Aut(\CY)$ is also closed in $\scrE$, according to \cite[Lemma 4.4.4]{Brion2}; thus, $Aut(\CY)$ is a union of connected components of $\scrE$. 
Note further that this method can not be implemented to treat algebraic stacks. This is because the construction of Hilb scheme for algebraic stacks is problematic. The reason is that the graph $\Gamma_f$ of a morphism $\CY\to \CY'$ of algebraic stacks is not closed in general, and requiring this is infact too much restrictive. For example one may observe that this assumption for the graph of diagonal $\Delta:\CY\to\CY\times_S\CY$ implies that $\CY$ is representable by an algebraic space. To handle the case of algebraic stacks one needs to make use of more sophisticated techniques such as deformation theory of 1-morphisms of algebraic stacks and etc.; e.g. see \cite{Ao}, \cite{Ols2} and \cite{HR}. 
\end{remark}

\subsection{Local Model For $\Sigma_\varsigma$}

Recall that for a stack $\CY$ over $S$, one defines the tangent bundle $\CT_\CY$ via the following functor
$$
T/S\mapsto \CY(T[\epsilon]),
$$
where $T[\epsilon]:=T\times_\BZ\BZ[\epsilon]$ with $\epsilon^2=0$.

One can see that when $\CY$ is an Artin stack then the same holds for $\CT_\CY$. The stack $\CT_\CY$ is equipped with the projection morphism $\CT_\CY\to \CY$ and zero section $\CY \to \CT_\CY$, that are induced by $\BZ[\epsilon]\to \BZ$, $\epsilon \mapsto 0$, and inclusion $\BZ\to\BZ[\epsilon]$ respectively.
Consider the groupoid $\scrE(\CT_\CY):=\ul\Hom(\CT_\CY,\CT_\CY)$ of endomorphisms of $\CT_\CY$. The projection morphism and zero section define a projection morphism $\scrE(\CT_\CY)\to\scrE$ and an obvious section $\scrE\to\scrE(\CT_\CY)$. Consider the following morphism
$$
d: \scrE\to \scrE(\CT_\CY)
$$
given by sending $\varsigma$ to its derivation $d\varsigma$. We let $\scrE_{d=0}$ denote the pull back $\scrE\times_{d,\scrE,0}\scrE(\CT_\CY)$.

\begin{theorem}\label{ThmLocModTFB}

Let $\ul\CH:=(\CH,char,pr^\leftarrow,pr^\rightarrow)$ be a $HR$-tuple, see \ref{Def_HR}. Let $\Sigma_\varsigma$ denote the fiber of $\Sigma(\ul\CH)$ over $\varsigma \in \scrE$. Let $z$ be a point in $\Sigma_\varsigma$, and set $y=pr^\leftarrow(z)$. Assume that

\begin{enumerate}
\item
$\CY$ is an Artin stack, smooth at $y$, which admits an \'etale neighborhood $\CU_y\to\CY$ at $y$,

\item
$y$ lies in the vanishing locus of $d\varsigma$, and

\item

there exist a family $\FF$ over $X$ and an \'etale neighborhood $\CU_{(y,x)}$ of $(y,x)$ that trivializes $\CH$ over $\CY\times X$ in the following sense

$$
\CH \times_{\CY \times X} \CU_{(y,x)}\tilde{\to} (\FF\times \CY)\times_{\CY\times X}\CU_{(y,x)}.
$$

\end{enumerate}  

\noindent
Then there is an \'etale neighborhood $U_z$ of $z$ and a roof of \'etale morphisms

\[ 
\xygraph{
!{<0cm,0cm>;<1cm,0cm>:<0cm,1cm>::}
!{(0,0) }*+{U_z}="a"
!{(-1.5,-1.5) }*+{\Sigma_\varsigma}="b"
!{(1.5,-1.5) }*+{\FF,}="c"
"a":^{\text{\'et}}"b" "a":_{\text{\'et}}"c"
}  
\]

\end{theorem}

\begin{proof}

Set $U=\CU_(y,x)$ and 

$$
U':=\CH|_U\cong \FF \times \CY|_U=:U''.
$$
\noindent
Consider the following diagram 

\[ \xygraph{
!{<0cm,0cm>;<1cm,0cm>:<0cm,1cm>::}
!{(.75,2)}*+{\CH}="a"
!{(.75,0)}*+{\CY\times\CY\times X}="n"
!{(3,-1)}*+{\CY\times X}="o"
!{(3,1)}*+{\Sigma_\varsigma}="p"
!{(-1,1)}*+{\wt U_z}="q"
!{(-2.25,2)}*+{U'}="b"
!{(.75,-2.25)}*+{\CY\times X}="c"
!{(-2.25,-2.25) }*+{U}="d"
!{(.75,-4)}*+{\FF\times \CY}="e"
!{(-2.25,-4) }*+{U''}="f"
!{(.75,-5.5) }*+{\FF}="h"
"b":^{\et}"a"
"d":^{\et}"c"
"f":^{\et}"e"
"p":"a"
"p":"o"
"o":_{\id\times\varsigma\times\id}"n"
"q":^{\et}"p"
"q":"b"
"a":"n"
"n":_{pr_2\times \id}"c"
"b":"d"
"e":_{pr_1}"h"
"e":"c"
"f":"d"
"b":@/_4em/_{\phi'}^{\cong}"f"
}  
\]

\noindent
We need to check that the composition of the maps
$$
\CD
\wt U_z @>>> U'@>\cong>> U'' @>\et>> \FF \times \CY @>pr_1>> \FF
\endCD
$$
is \'etale. Note that we may realize $\wt U_z$ as a substack of $U'$ given by
$$
\wt U_z:=\equi( \varsigma \circ f, g\colon  U' \rightrightarrows \CY).
$$
 Here $f$ is the morphism induced by the projection $pr^\leftarrow: \CH\to \CY$ and $g$ is given by $U'\cong U''\to \FF\times \CY$ followed by the projection to the second factor.\\
Let $v\in \wt U_z$ be a point and let $w\in U'$ and $y=\varsigma\circ f(w)=g(w)\in \CY$, as well as $t\in \FF$ be its images. By smoothness of $\CY$, we may further reduce to the case that $\CY$ is an affine space. Namely take an affine open neighborhood $Y'$ of $y$ in $\CY$ which admits an \'etale morphism $\pi\colon Y'\to\wt Y$ to some affine space $\wt Y=\BA^m_S=\Spec\CO_S [z_1,\ldots,z_m]$, and consider an affine neighborhood $T'$ of $t$ which we write as a closed subscheme of some $\wt T=\BA^\ell_{\CO_S}$. Replace $Y$ by the affine neighborhood $Y'$ of $y$ and $U'$ by an affine neighborhood $W'$ contained in $(\varsigma\circ f)^{-1}(Y')\cap g^{-1}(Y')\cap \iota^{-1}(T'\times Y')$. Then $\wt U_z':=W'\times_{U'} \wt U_z=W'\times_{(\sigma_Zf,g),Y'\times Y',\Delta}Y'$ is an open neighborhood of $v$ in $\wt U_z$. We may extend the \'etale morphism $\iota\colon W'\to T'\times Y'$ to an \'etale morphism $\tilde\iota\colon\wt W\to\wt T\times Y'$ with $\wt W\times_{\wt T}T'=W'$. We also extend $\pi\circ f\colon W'\to\wt Y$ to a morphism $\tilde f\colon\wt W\to \wt Y$ and $g$ to $\tilde g \colon\wt W\to \wt Y$, and set $\wt{\wt U}:=\equi(\varsigma \circ\tilde f,\tilde g\colon\wt W\rightrightarrows\wt Z)=\wt W\times_{(\sigma_{\wt Z}\tilde f,\tilde g),\wt Z\times\wt Z,\Delta}\wt Z$. Since $\Delta\colon Y'\to Y'\times_{\wt Y} Y'$ is an open immersion, also the natural morphism 
\[
\wt U\;\longto\; \wt{\wt U}\times_{\wt T}T'\;=\;W'\times_{(\varsigma\tilde f,\tilde g),\wt Y\times\wt Y,\Delta}\wt Y\;=\;W'\times_{(\varsigma f,g),Y'\times Y',\Delta}(Y'\times_{\wt Y}Y')
\]
is an open immersion. Since $\wt W$ is smooth over $\wt T$ of relative dimension $m$ and $\wt{\wt U}$ is given by $m$ Artin-Schreier type equations 
$$
\tilde g^*(z_j)- \varsigma^\ast \circ \tilde f ^* (z_j),
$$ 
and $y$ lies in the vanishing locus of $d\varsigma$ and thus the above equations have linearly independent differentials $d\tilde g^*(z_j)$, we observe that $\wt{\wt U}\to\wt T$ is \'etale  according to the Jacobi-criterion \cite[\S 2.2, Proposition~7]{BLR}.

\end{proof}


The above theorem has the following immediate corollaries that partly describe global and local geometry of $\Sigma_\varsigma$ for $\varsigma$ in $\scrE_{d=0}$.\\

\noindent
Let $S=\Spec k$. Let $\ul\CH$ be a $HR$-tuple. Assume further that there is $y_0$ in $\CY(k)$ and that $\CH\to\CY\times X$ is a fiber bundle with fiber $\FF:=\CH_{y_0}$ for the \'etale topology on $\CY\times_S X$. 
Suppose $\CH$ admits a stratification $\{\CH_\lambda\}_\lambda$, that further induces stratification $\{(\Sigma_\varsigma)_\lambda\}_\lambda$ of $\Sigma_\varsigma$. One can easily observe that

\begin{corollary}\label{CorICSheafNablaH}
Keep the notation and assumptions in Theorem \ref{ThmLocModTFB} together with the above notation and assumptions. Let $\varsigma$ be in $\scrE_{d=0}$. The IC-sheaf $IC(\Sigma_\varsigma)$ is the restriction of $IC(\CH)$ up to some shift and Tate twist.

\end{corollary}

\begin{proof} Let $\Sigma_\ul c$ denote the fiber above the constant morphism $\ul c: \CY\to\CY$, sending $\CY$ to the point $y_0$. According to \ref{ThmLocModTFB}, we may replace $\FF$ by $\Sigma_\ul c$. Therefore we may assume that $\FF$ is equipped with a natural stratification induced by $\{\CH_\lambda\}_\lambda$ of $\CH$. Consider the following diagram
\[ \xygraph{
!{<0cm,0cm>;<1cm,0cm>:<0cm,1cm>::}
!{(1.5,2)}*+{\CH}="a"
!{(1.5,0)}*+{\CY\times\CY\times X}="n"
!{(4,-2.25)}*+{\CU;}="o"
!{(4,2)}*+{\CH|_\CU}="p"
!{(5.5,2)}*+{\FF\times\CY|_\CU}="r"
!{(8,2)}*+{\FF}="t"
!{(4.5,2)}*+{\cong}="u"
!{(-2.25,2)}*+{\Sigma_\varsigma}="b"
!{(1.5,-2.25)}*+{\CY\times X}="c"
!{(-2.25,0) }*+{\CY\times X}="d"
!{(-4.25,2)}*+{\CU_z}="e"
"e":"b"
"b":"d"
"b":"a"
"p":^{i'}"a"
"r":^{pr_1}"t"
"p":"o"
"o":"c"
"a":"n"
"n":"c"
"d":^{(id, \varsigma)\times id}"n"
"e":@/^3em/^{i}"p"
}  
\]
see the proof of the above theorem \ref{ThmLocModTFB}.
The stratification on $\Sigma_\varsigma$ is induced by that of $\CH$. Moreover by  Theorem \ref{ThmLocModTFB} the smooth open stratum $\Sigma_\varsigma^\circ$ lies inside the pull back of the open smooth stratum  $\CH^\circ$ of $\CH$, the statement is obvious over the open stratum. By the Theorem \ref{ThmLocModTFB} we have $IC(\CU_z)=(pr_1\circ i)^\ast IC(\FF)$ which equals $i^\ast IC(\CH|_\CU)$ up to shift and Tate twist by $\dim \CY$. By \'etaleness of $i'$ the later coincides $i^\ast i'^\ast IC(\CH)$ which equals the restriction of $IC(\CH)$ to ${\CU_z}$.

\end{proof}

\begin{corollary}\label{CorSingularities}
Keep the notation and assumptions in \ref{CorICSheafNablaH}. The stalks $\wh\CO_{\Sigma_\varsigma,z}$ at the point $z$ satisfy Serre's condition $S_i$ (resp. $R_i$) if the points of the fiber $\Sigma_\ul c$ satisfy $S_i$ (resp. $R_i$).

\end{corollary}

\subsection{Lifting of Fibers}\label{SubsecLifing}

Let $S$ be the formal spectrum of a complete discrete valuation ring $R$, with special point $s$ and generic point $\eta$. Let $k=\kappa(s)$ (resp. $\kappa(\eta)$) denote the residue fields at $s$ (resp. $\eta$). Set $\ol S:=\Spec k$. Choose a separable closure $\ol\eta$ of $\eta$, and let $\ol s$ be the residue field of the normalization $\tilde{S}$ of $S$ in $\kappa(\eta)$.

\begin{lemma}\label{LemLiftOfCGExist}
Let $\CW\to \ul\Hom$ be a $S$-morphism from a smooth stack $\CW$.
\forget{Let $\{\varsigma\}_n$ be a compatible set of endomorphsm of $\CY\times_R R/\Fm_R^n$.} Let $\ol\varsigma$ be a $\ol S$-point of $\ul\Hom$ which comes from a $\ol S$-point of $\CW$. Then there is a stack $\wh \Sigma_\varsigma$ over $\Spf R$ which lifts $\Sigma_{\ol\varsigma}$ over $\Spf R$. Assuming further that $\CH$ and $\CY$ are proper, then there is a $Gal(\ol\eta/\eta)$-equivariant isomorphism

$$
\Koh^i(\wh\Sigma_\ol\eta,\ol\BQ_\ell) \cong \Koh^i(\Sigma_{\ol\varsigma},R\Psi(\ol\BQ_\ell)).
$$
Here $R\Psi(\ol\BQ_\ell)$ denotes the corresponding sheaf of nearby cycles.
\end{lemma}

\begin{proof}
The point $\ol\varsigma: \ol S \to \ul\Hom$ comes from $\ol S$-point of the smooth stack $\CW$. By infinitesimal criterion of smoothness it gives a compatible set of morphisms $\varsigma_n:\CY_n:=\CY\times_R R/\Fm_R^n\to \CY_n$\forget{which further gives $\varsigma: \Spf R\to \wh\scrE$}. This consequently yields a compatible system $ \Sigma_{\varsigma_n}\hookrightarrow \Sigma_{\varsigma_{n+1}}$ of closed immersions of algebraic stacks, which accordingly define a formal algebraic stack $\wh\Sigma_\varsigma$ over $\Spf R$. Let $\wh\CY:=\CY\hat{\times}\Spf R$ and $\wh\CH:=\wh\CY\hat{\times}\Spf R$. Assuming that $\CH$ and $\CY$ are proper we see that $\wh \CH\to\wh\CY\times \wh\CY$ and also $\Sigma_\ol\varsigma$ are proper, thus we see by Grothendieck existence theorem that $\wh \Sigma_\varsigma$ lifts to a pojective family over $\Spec R$. The isomorphism 
$$
\Koh^i(\wh\Sigma_\ol\eta,\ol\BQ_\ell) \cong \Koh^i(\Sigma_{\ol\varsigma},R\Psi(\ol\BQ_\ell))
$$
of the cohomology groups follows from basic properties of the sheaf of nearby cycles; see for example \cite[Theorem 10.1]{AEK}. For the Grothendiek existence theorem in the context of algebraic stacks see \cite[Theorem 1.4]{Ols1}
\end{proof}


\section{The family $\Sigma\to\scrE$ arising from Hecke stack}\label{SectShtukasManifestations}

In this section we focus on particular examples of the construction we discussed in Section \ref{SubsectTwoFiberBundles}. Before doing this, we need to recall some further preliminary materials.

\subsection{The moduli stack $\scrH^1(C,\FG)$}

 Let $X\to S$ be a projective flat morphism of schemes. Notice that we later restrict ourselves to the case that $X$ is a (relative) curve $C$ (over $S$). Let $\scrH^1(X,\FG)$ be the stack classifying $\FG$-bundles on $X$.\\
Assume that $\FG$ admits a representation $\iota: \FG\to \GL(\CV_0)$, where $\CV_0$ is a vector bundle over $X$, and such that it fulfills the following requirment

\begin{eqnarray}\label{EqCondQuotient}
\text{there is a scheme $Y$ affine and of finite type over $X$ with an action $\FG\times_X Y\to Y$}\\ \text{ of $\FG$ and a $\GL(\CV_0)$-equivariant open immersion $\GL(\CV_0)/\FG\into Y$.}\nonumber
\end{eqnarray}

\noindent
To see up to what extent the above condition can be served see \cite{AH_Global}. 

The representation $\iota$ induces a morphism $\iota: \scrH^1(X,\FG)\to \scrH^1(X,\GL(\CV_0))$ of $S$-stacks.  For a scheme $T$ and a morphism $T\to \scrH^1(X,\GL(\CV_0))$, corresponding to $\GL(\CV_0)$-torsor $\CG$, one forms the following 2-Cartesian diagram

$$
\CD
\pi_\ast(\CG/\FG)@>>>T\\
@VVV @VV{\CG}V\\
\scrH^1(X,\FG)@>>> \scrH^1(X,\GL(\CV_0))
\endCD
$$
  
The above condition ensures that $\pi_\ast(\CG/\FG)$ is a quasi-affine $S$-scheme of finite presentation; see \cite[Theorem~2.6]{AH_Global}. This way one reduces the study of $\scrH^1(X,\FG)$ to the well-known case where $\FG=\GL(\CV_0)$.

Let us state the following basic result. 

\begin{theorem}\label{ThmBunG}

 Let $X\to S$ be a projective flat morphism of schemes. Let $\FG$ be as above. Then $\scrH^1(X,\FG)$ is an algebraic $S$-stack locally of finite presentation.

\end{theorem}

\begin{proof}
The theorem is well known when $\FG$ is constant, for a split reductive group $G$, and $X=C$, where $C$ is a smooth projective curve over a perfect field $k$. When $\FG$ is a parahoric group scheme over a smooth projective curve $C$ over a perfect field $k$ see Heinloth \cite[Proposition 1]{Heinloth}. For more general case where $\FG$ is a flat affine group scheme of finite type over $C$, a proof is given in \cite[Theorem 2.5]{AH_Global}. The idea is using the method discussed above and showing that a flat affine group scheme of finite type over the curve $C$ satisfies the above condition \ref{EqCondQuotient}. The statement for the general relative case is similar, except that we do not have the relative version of the \cite[Proposition~2.2]{AH_Global} which ensures the existence of the representation $\iota:\FG\to\GL(\CV_0)$ that satisfies the above condition \ref{EqCondQuotient}. Note however that this is obvious for the constant reductive case.   Also for the relative case see \cite[Theorem~1.0.1]{Wang}.
\end{proof}

\begin{remark}\label{Rem_StackGBun}

When $X=C$, where $C$ is a smooth projective family of curves with geometrically reduced, connected fibres over $S$, and $\FG$ is smooth, then the stack $\scrH^1(C,\FG)$ admits an open covering $\{\CU_\alpha\}_{\alpha\in I}$ by smooth algebraic substacks of finite presentation over $S$; e.g. see \cite[Theorem~1.0.1]{Wang} or \cite[Theorem~2.5]{AH_Global}.
Furthermore it's diagonal morphism is schematic, affine and of finite presentation.
  
\end{remark}

\subsubsection{Heinloth-Schmitt $a$-stability condition}\label{SubsectionHeinlothSchmitt}

Note that in \cite{HeiSch} the authors establish this theory for the case that $S=\Spec k$ for a finite field $k$ and $\FG$ is constant, i.e. $\FG=G\times_{\BF_q}C$ for a split reductive group $G$ over $k$. Then they explain that their theory carries over to the case where $G$ is a reductive group over an integral ring $R$, finitely generated over $\BZ$, up to some modifications, see \cite[Remark 3.2.4]{HeiSch}. This is essential for the techniques they implement in their article, see proof of \cite[Corollary 3.3.4 and Theorem 3.3.5]{HeiSch}. We don't know up to what extent they may remain valid for more general $\FG$. According to this, whenever we make use of the stability condition, we implicitly assume that $\FG$ is constant for a reductive  group $G$ over $S$. \\

Here for the convenience of the reader we briefly recall the $\ul a$-(semi)stability result of Heinloth-Schmitt from \cite{HeiSch}. 

\begin{definition-remark}
\begin{itemize}
\item[a)] To recall the $\ul a$-stability condition we recall the definition of the $G$-bundles with flagging of type $(\ul x,\ul P)$.
Here $\ul x = (x_i)_{i=1,...,b}$ denote a finite set of distinct $k$-rational points of $C$, and $\ul P = (P_i)_{i=1,...,b}$ denotes a tuple of parabolic subgroups of $G$. A principal $G$-bundle with a flagging of type $(\ul x,\ul P)$ is a tuple $(\CP,\ul s)$ that consists of a principal $G$-bundle $\CP$ on $C$ and a tuple $\ul s=(s_1,...,s_b)$ of sections $s_i: {x_i}\to(\CP\times_C \{x_i\})/P_i$. The category of principal $G$-bundles with a flagging of type $(\ul x,\ul P)$ form the smooth algebraic stack $\scrH^1(C,\FG,\ul x,\ul P)$; see \cite[Lemma~3.2.2]{
HeiSch}. \\

\item[b)] For an algebraic group $P$
let $X^\ast(P)$ denote the corresponding group of characters. Let $\check{X}^\ast(P)_{\BQ,+}$ denote the set of all elements $a\in \check{X}^\ast(P)_\BQ$  such that $a(\det_{P'}\otimes \det_{P})< 0$ for all $P\subseteq P'$. \\
Let $(\CP,\ul s)$ be a flagged principal $G$-bundle and $\CP_{x_i,P_i}$
the $P_i$-torsor over $x_i$ defined by $s_i$, $i = 1,\dots,b$. Set $P_{s_i}:= Aut_{P_i}(\CP_{x_i,P_i}) \subseteq \Aut_G(\CP_{x_i})$ the corresponding parabolic subgroup. There are canonical isomorphisms $X^\ast(P_i)_\BQ\cong X^\ast(P_{s_i})_\BQ$ and $\check{X}^\ast(P_i)_{\BQ,+}= \check{X}^\ast(P_{s_i})_{\BQ,+}$, $i = 1,\dots ,b$. For $a_i \in \check{X}^\ast (P_i)_{\BQ,+}$ we let $a_{s_i}$
denote the corresponding element in $\check{X}^\ast(P_{s_i})_{\BQ,+}$. For a parabolic subgroup $Q$ of $G$, a character $\chi$ of $Q$, and a reduction $\CP_Q$ of $\CP$ to $Q$, we get in each point $x_i$ a parabolic subgroup $Q_i$ in $\Aut(\CP_{x_i})$ and a character $\chi_{s_i}$ of that parabolic subgroup, $i=1,\dots ,b$. Fix $\ul a \in \prod_i X^\ast(P_i)$. For a parabolic subgroup $Q$ of $G$ and a reduction $\CP_Q$ of $\CP$ to $Q$, define

\begin{eqnarray*}
\ul a-deg(\CP_Q) \colon  X^\ast(Q) &\to&  \BQ\\
\chi &\mapsto & \deg(\CP_Q(\chi))+\sum_i\langle \chi_{s_i},a_{s_i}\rangle.\\  
\end{eqnarray*}

Here  $\CP_Q(\chi)$ is the line bundle on $C$ given by pushing forward the principal $Q$-bundle $\CP_\CQ$ via the character $\chi$, and $\chi_{s_i}$ is a character of $Q_i \subseteq \Aut(P_{x_i})$ in each point $x_i$.  To see that the pairing $\langle \chi_{s_i},a_{s_i}\rangle$ is well-defined see \cite[Remark~4.1.2 iii)]{HeiSch}.

\item[c)] A flagged principal $G$-bundle $(\CP,\ul s)$ is called $\ul a$-semistable (resp. stable), if for any parabolic subgroup $Q \subseteq G$
and any reduction $\CP_Q$ of $\CP$ to $Q$, the following holds
$$
\ul a-\deg(\CP_Q)\leq 0
$$ 
(resp. $\ul a-\deg(\CP_Q)< 0$).

Denote by $\scrH^1(C,\FG,\ul x, \ul P)^{\ul a-(s)s}$, the substack of the moduli stack $\scrH^1(C,\FG,\ul x, \ul P)$ parametrizing $\ul a$-(semi)stable flagged principal $\FG$-bundles of type $(\ul x,\ul P)$. We further use the notation $\scrH^1(C,\FG)_{(\ul x, \ul P)}^{\ul a-(s)s}$ to denote the (scheme theoretic) image of $\scrH^1(C,\FG,\ul x, \ul P)^{\ul a-(s)s}$ under the projection 
$$
\scrH^1(C,\FG,\ul x, \ul P)^{\ul a-(s)s}\to\scrH^1(C,\FG);
$$
see \cite[Tag 0CMH, Lemma 98.37.3.]{Stacks} for existence of the scheme theoretic image. Note that for the scheme theoretic image $Z$ of a quasi-compact morphism of algebraic stacks $f:\CX\to\CY$, one can observe that $|Z|$ is the closure of the image of $|f|$; see \cite[Tag 0CMH, Lemma 98.37.6.]{Stacks}.

\end{itemize}
\end{definition-remark}

Let us now state the following important result of Heinloth and Schmitt which was the key point in their proof of the purity of $\Koh^\ast(\scrH^1(C,\FG),\BQ_\ell)$.

\begin{theorem}\label{ThmHeinlothSS}
For any substack $\CU$ of $\scrH^1(C,\FG)$, and any integer $i>0$, there is a type $(\ul x,\ul P)$ and stability condition $\ul a$ such that $\scrH^1(C,\FG,\ul x, \ul P)^{\ul a-(s)s}$ is smooth proper, and it's image contains $\CU$ (in particular the closure $\ol \CU$ in $\scrH^1(C,\FG)$ is proper). Furthermore the codimension of $\scrH^1(C,\FG,\ul x, \ul P)^{\ul a-(s)s}$ in $\scrH^1(C,\FG)_{\ul x, \ul P}$ is $>i$. 
\end{theorem}

\begin{proof}
This is \cite[Theorem 3.2.3]{HeiSch}.
\end{proof}

\begin{proposition}\label{PropAoki}
Fix open substacks $\CU, \CV\subseteq \scrH^1(C,\FG)$ of finite type. The stack 
$$
\Hom(\scrH^1(C,\FG)_{(\ul x, \ul P)}^{\ul a-(s)s},\scrH^1(C,\FG))
$$ (resp. $\scrE_{(\ul x, \ul P)}^{\ul a-ss}:= Hom(\scrH^1(C,\FG)_{(\ul x, \ul P)}^{\ul a-(s)s},\scrH^1(C,\FG)_{(\ul x, \ul P)}^{\ul a-(s)s})$, resp. $\Hom(\ol\CU,\ol\CV)$) is an algebraic stack for appropriate choice of $(\ul x, \ul P)$. Furthermore, $\Hom(U,V)$, for affine charts $U\to \scrH^1(C,\FG)$ and $V\to \scrH^1(C,\FG)$ (resp. $U\to \scrH^1(C,\FG)_{(\ul x, \ul P)}^{\ul a-(s)s}$ and $V\to \scrH^1(C,\FG)_{(\ul x, \ul P)}^{\ul a-(s)s}$), is smooth.
\end{proposition}

\begin{proof}

The algebraicity of these stacks follow from the above theorem \ref{ThmHeinlothSS} and the main result of \cite{HR}; see also \cite{Ols2}. The smoothness follows from infinitesimal criterion for smoothness and the fact that $\scrH^1(C,\FG,\ul x,\ul P)$ is a locally trivial bundle over the smooth Artin stack $\scrH^1(C,\FG)$ whose fibers are isomorphic to $\prod_{i=1}^b(G/P_i)$. 

\end{proof}

\subsection{Construction of $\Sigma\to \scrE$ for $\CH=Hecke_n(C,\FG)$}

In this subsection we first recall the definition of the Hecke stack. We view it as a two fiber bundle over $\scrH^1(C,\FG)$. This yields a family $\Sigma\to \scrE$ via the construction explained in section \ref{SectGeneralConstruction}. We then discuss certain fibers of this family and we recall the notion of boundedness condition from \cite{AH_Local}, \cite{AH_Global} and \cite{AH_LM}. We further observe that a boundedness condition gives rise to certain Hecke classes in Chow group; see Proposition \ref{PropHeckOverBunG}. We finally construct a local model roof for the family $\Sigma\to \scrE$, see Proposition \ref{Prop_LocalModel}.

\begin{definition}\label{Hecke}
For a natural number $n$, let $Hecke_n(C/S,\FG)$ be the stack fibered in groupoids over the category of $S$-schemes, whose category of $T$-valued points consists of tuples $\bigl(\CG,\CG',\ul\charsect,\tauGlob\bigr)$, where
\begin{itemize}
\item[--] $\CG$ and $\CG'$ are in $\scrH^1(C,\FG)(T)$,
\item[--] $\ul\charsect:=(\charsect_i)_i \in C^n(T)$ is an n-tuple of sections, and
\item[--] $\tau\colon  \CG_{|_{{C_T}\setminus{\Gamma_{\ul\charsect}}}}'\isoto \CG_{|_{{C_T}\setminus{\Gamma_{\ul\charsect}}}}$ is an isomorphism.
\end{itemize}
Forgetting the isomorphism $\tau$ defines a morphism 
\begin{equation}\label{EqFactors}
Hecke_n(C/S,\FG)\to \scrH^1(C,\FG)\times_S \scrH^1(C,\FG) \times_S C^n.
\end{equation}
\noindent

\end{definition}

\begin{proposition}\label{PropHeckeisInd_Quasi_Proj}
Let $\FG$ be a flat affine group scheme of finite type over $C$. Furthermore assume that it admits a faithful representation $\rho:\FG\to\SL(\CV_0)$, for a vector bundle $\CV_0$ over $C$, with affine (resp. quasi-affine) quotient $\SL(\CV_0)/\FG$. The stack $Hecke_n(C/S,\FG)$ is ind-algebraic stack, ind-projective (resp. ind-quasi-projective) over $C^n\times \scrH^1(C,\FG)$; see \cite[Propostion 3.9]{AH_Global}. Furthermore it is ind-projective if $\FG$ is parahoric.\\
\end{proposition}
\begin{proof}
See \cite[Proposition~3.9]{AH_Global}. Note however that the argument given in loc. cit. is stated for $S=\Spec k$, but nevertheless, one can literally follow the same lines to prove the more general statement.

\end{proof}

\noindent
\begin{notation-remark}\label{Not_Rem1}
\begin{enumerate}

\item\label{Rem-Not a} Consider the $HR$-tuple $\ul\CH_n(C,\FG):=(\CH_n(C,\FG), char, pr^\leftarrow,pr^\rightarrow)$, where\\[1mm]
\noindent
\begin{enumerate}

\item[--] $\CH_n(C,\FG):=Hecke_n(C/S,\FG)$,
\item[--] $char: \CH_n(C,\FG)\to C^n$ $(\CG,\CG',\ul s, \tau)\mapsto \ul s$,
\item[--] $\CY(C,\FG):=\scrH^1(C,\FG)$,
\item[--] $pr^{\leftarrow}: \CH_n(C,\FG) \to \CY(C,\FG)$ (resp. $pr^{\rightarrow}: \CH_n(C,\FG) \to \CY(C,\FG)$) the projection to the first (resp. second) factor, see \ref{EqFactors},

\end{enumerate}

For $HR$-tuple $\ul\CH:=\ul\CH_n(C,\FG):=(\CH_n(C,\FG), char, pr^\leftarrow, pr^\rightarrow)$ as above, we set $\Sigma_n(C,\FG):=\Sigma(\ul\CH)$; see Definition \ref{Def_HR}. \\

\noindent
When the curve $C$ and the group $\FG$ are obvious from the context, we remove $(C,\FG)$ from our notation. Also when $n$ is clear from the context we write $\CH=\CH_n$ and $\Sigma=\Sigma_n$.

\item\label{Rem-Not b}
We view $\Sigma=\Sigma(C,\FG)$ as a family $\Sigma\to \scrE :=\scrE(C,\FG):=\ul\Hom(\CY,\CY)$.
 and we denote by $\Sigma_\varsigma$, the fiber above $\varsigma\in\scrE$.  
We call an object $\ul\CG$ in $\Sigma_\varsigma(S)$ \emph{a (global) $\varsigma$-$\FG$-shtuka} over $S$.

\item\label{Rem-Not c}
For a relative Cartier divisor $D$ on $C$, we denote by $\CY_D:=\scrH_D^1(C,\FG)$ the stack over $\scrH^1(C,\FG)$, whose $T$-points classifies $\FG$-bundles on $C_T$ together with $D$-level structures, i.e. a trivialization $\CG|_{D\times T}\to \FG\times_{C}C_T|_{D\times T}$. Using this one can further equip the other  constructions with level $D$-structure in an obvious way. We use the subscript $D$ in our notation $\CH_D$, $\Sigma_D$ and etc. to illustrate that the objects parameterized by corresponding moduli stacks are equipped with D-level structures.

\item\label{Rem-Not d} Assume that $\FG=G\times_S C$ for split reductive group $G$ over $S$. Fix a type $(\ul x, \ul P)$; see subsection \ref{SubsectionHeinlothSchmitt}. We similarly use the notation $\CH_{(\ul x, \ul P)}^{\ul a-(s)s}$ for the restriction of $\CH$ to $\CY_{(\ul x, \ul P)}^{\ul a-(s)s}:=\scrH^1(C,\FG)_{(\ul x, \ul P)}^{\ul a-(s)s}$ under $pr^\leftarrow: \CH\to\CY$. Regarding the procedure described in subsection \ref{SubsectTwoFiberBundles}, we obtain the family $\Sigma_{(\ul x, \ul P)}^{\ul a-(s)s}$ over 
$$
\ul{\Hom}(\scrH^1(C,\FG)_{(\ul x, \ul P)}^{\ul a-(s)s}, \scrH^1(C,\FG)).
$$
Note that the latter stack is an algebraic stack for appropriate choice of $(\ul x, \ul P)$, see Proposition \ref{PropAoki}.

\item\label{Rem-Not e} One may alternatively require that the $\FG$-bundles occurring in the above parts are semistable in the sense of \cite{BaPa}. For this one requires $G$ to be semi-simple. We set $\CH^{ss}=Hecke_n^{ss}$, $\CY^{ss}$, $\scrE^{ss}$ and etc. for the corresponding moduli spaces. 
Note that the moduli space  $\CY^{ss}$ is projective, according to \cite{BaPa} and \cite{BaSe}. Consequently, $\scrE^{ss}$ is representable by a scheme locally of finite type. See Remark \ref{Remark_HilbStack} and also \cite{Brion1} for some details on the structure of the endomorphisms of projective varieties.

\end{enumerate}

\end{notation-remark}

\subsection{Specific fibers of $\Sigma\to \scrE$}\label{Subsec_SpecificFibers}

 As we will see bellow, some interesting moduli spaces appear in the fibers of the family $\Sigma\to\scrE$.

\noindent
\subsubsection{Global affine Grassmannian, Boundedness Conditions and Local Model}\label{SubsectBDGrassmannian}

\begin{definition}\label{DefBDGr}
Fix a $\FG$-bundle $\CG_0$ in $\CY:=\scrH^1(C,\FG)(S)$. Take $\varsigma$ to be the  constant morphism  
$$
\CD
\ul\Theta_{\CG_0}:\CY\to\CY,\\ 
~~~~~y\mapsto \CG_0. 
\endCD
$$

We denote by $\Sigma_{\ul\Theta_{\CG_0}}$
the fiber of the family $\Sigma\to\scrE$. The stack $\CY$ has an especial $S$-point corresponding to the trivial $\FG$-bundle, we denote the corresponding constant morphism by $\Theta$. The fiber $\Sigma_{\ul\Theta}$ is called the \emph{(relative) Beilinson-Drinfeld} affine Grassmannian. We also use the notation $GR_n(C,\FG)$ for $\Sigma_{\ul\Theta}$.\\

\end{definition}

\begin{proposition}\label{Prop_Gr_nisindQProj}
Let $\FG$ be a flat affine group scheme of finite type over $C$. Furthermore assume that it admits a faithful representation $\rho:\FG\to\SL(\CV_0)$, for a vector bundle $\CV_0$ over $C$, with affine (resp. quasi-affine) quotient $\SL(\CV_0)/\FG$. Then the fiber $\Sigma_{\ul\Theta_{\CG_0}}$ is an ind-scheme ind-projective (resp. ind-quasi-projective) over $C^n$. In particular when $\FG$ is parahoric the Beilinson-Drinfeld affine Grassmannian $GR_n(C,\FG)$ is an ind-scheme ind-projective (resp. ind-quasi-projective) over $C^n$. 
\end{proposition}

\begin{proof}

This follows from Proposition \ref{PropHeckeisInd_Quasi_Proj}.
\end{proof}

\begin{definition-remark}
The (relative) loop group $\CL_n\FG$ (resp. positive loop group $\CL_n^+\FG$) is the space corresponding to the following functor
$$
T\mapsto\{(s_1,\dots,s_n,\alpha); s_i\in C(T) ~and~ \alpha\in \FG(\dot{\BD}(\Gamma_{\ul s})\}
$$
(resp. $$
T\mapsto\{(s_1,\dots,s_n,\alpha); s_i\in C(T)~and~\alpha\in \FG(\BD(\Gamma_{\ul s})\})
$$

According to gluing lemma of Beauville-Laszlo, \cite{B-L}, one observes that $\CL_n\FG$ and $\CL_n^+\FG$ operate on $\Sigma_\Theta$ via changing the trivialization.
 
\end{definition-remark}

Note that one can use Beilinson-Drinfeld affine Grassmannian  $GR_n(C,\FG)$ to locally trivialize the family $\CH\to C^n\times\CY$.

\begin{proposition}\label{PropHeckeLocMod}
Consider the stacks $\CH$ and $GR_n(C,\FG) \times \CY$ as families over $C^n \times \CY$, via the projections $(\CG,\CG',\ul s,\tau) \mapsto (\ul s,\CG')$ and $(\CG,\ul s,\tau) \times \CG' \mapsto (\ul s,\CG')$ respectively. They are locally isomorphic with respect to the \'etale topology on $C^n \times \CY$. 

\end{proposition}

\begin{proof}
For the proof in the case $S=\Spec k$ see \cite[Proposition~2.0.11]{AH_LM} and \cite[Lemma~4.1]{Var}. The proof for the relative situation over a base scheme $S$ is similar.
\end{proof}

As we mentioned before, the moduli stack $Hecke_n(C,\FG)$ is  an ind-algebraic stack locally of ind-finite type. To provide a moduli stack which is (locally) of finite type, one may proceed by introducing boundedness conditions. There are various methods to establish such conditions. For a split reductive group $G$, Varshavsky uses an $n$-tuple $\ul\mu:=(\mu_i)$ of dominant coweights of $G$ to control the relative position of $\CG'$ and $\CG$ through $\tau:\CG'\to\CG$; see \cite[Definition 2.4]{Var}. Here we briefly recall the boundedness conditions that have been considered in \cite{AH_Global} and \cite{AH_LM} for a flat group scheme $\FG$ over $C$.

\begin{definition-remark}\label{Def_Rem_BC}

\begin{enumerate}
\item
Fix a faithful representation $\varrho: \FG\to\SL(\CV_0)\subseteq \GL(\CV_0)$ for some vector bundle $\CV_0$ of rank $r$ with quasi-affine quotient. Consider the induced morphism of stacks:\\
$$
\CD
\scrH^1(C,\FG)@>{\varrho_\ast}>>\scrH^1(C,\SL(\CV_0))@>\CV(-)>>\scrH^1(C,\GL(\CV_0))\cong Vect_C^r.
\endCD
$$
Let $\ul\omega:=(\omega_i)$ be an n-tuple of coweights of $\SL_r$ given as 
$$
\omega_i: x\mapsto \diag(x^{\omega_{i,1}},\dots,x^{\omega_{i,r}}),
$$
for integers $\omega_{i,1}\geq\dots\geq \omega_{i,r}$ with $\sum_\ell \omega_{i,\ell}=0$.\\

We say that a morphism $\tau:\CG\to \CG'$ between $\FG$-bundles $\CG$ and $\CG'$ over $C_T$, defined outside graph of the sections $s_i:T\to C$, is bounded by $\ul\omega$ if

$$
\wedge_{C_T}^\ell \tau^{-1} (\CV(\varrho_\ast\CG))\subseteq \wedge_{C_T}^\ell \CV(\varrho_\ast\CG')(\sum_{j=1}^\ell \omega_{i,j})\Gamma_{s_i},
$$
with equality when $\ell=r$. We denote by $Hecke^\ul\omega(C,\FG)$ the corresponding stack obtained by imposing the above boundedness condition. This further induces a boundedness condition on $\Sigma$. Similarly we use the notation $\Sigma^\ul\omega$ for the resulting bounded moduli stack.

\item
The above boundedness condition is not intrinsic as it depends to the choice of a representation. To provide an intrinsic boundedness condition, in \cite{AH_Global} and \cite{AH_LM}, we discussed another method. Namely, according to this method, a boundedness condition is given by a class of closed $\CL_n^+\FG$-stable subschemes $\CZ\subseteq \Sigma_\Theta\times_{C^n}\wt C$, where $\wt C$ is a smooth projective curve which is finite over $C$. Such a class of subschemes determine a minimal curve $C_\CZ$, which is called \emph{reflex curve}, over which the bounded moduli stack is defined; see \cite[Definition 3.1.3 and 4.3.2]{AH_LM}. To avoid the complications arising in the general set up, we assume that the reflex curve is $C$ itself, and the boundedness condition is given by $\CL_n^+\FG$-stable closed subschemes $\CZ\subseteq GR_n$. We in addition assume that $\CZ$ is proper and flat over $C^n$. We say that a morphism $\CG'\to\CG$ defined over $C_T\setminus \Gamma_\ul s$ is bounded by $\CZ$, if for every trivialization of $\CG$, the induced morphism $T\to GR_n$ factors through $\CZ$. This gives boundedness condition $\CZ$ on $Hecke_n(C/S,\FG)$, which further induces boundedness condition $\CZ$ on the moduli stack $\Sigma$. We denote the corresponding moduli spaces (stacks) obtained by imposing the boundedness condition $\CZ$, by $\CH_{n,D}^\CZ=Hecke_{n,D}^\CZ$, $\Sigma_D^\CZ$ and etc.

\end{enumerate}

\end{definition-remark}

\begin{remark}
When $\FG$ is parahoric then $GR_n(C,\FG)$ is ind-proper. In particular a closed subscheme $\CZ$ of $GR_n(C,\FG)$ is automatically projective; see Proposition \ref{Prop_Gr_nisindQProj}. But still flatness of $\CZ$ over $C^n$ is not obvious.
\end{remark}

\begin{proposition}\label{PropHeckOverBunG}
Assume that $\FG$ is smooth and $S=\Spec k$.
Let $\CZ$ be as in Definition-Remark \ref{Def_Rem_BC}b). Then the tuple $$\ul\CH_D^\CZ:=(\CH_{n,D}^\CZ, char, pr^\leftarrow, pr^\rightarrow)$$ is a $CHR$-tuple. Furthermore it induces a cycle $\phi_{D,\alpha}^\CZ$ in $Ch_{d+d'}(\scrH_D^1(C,\FG)_\alpha \times_S \scrH_D^1(C,\FG)_\alpha\times_S (C\setminus D)^n)$ for $D$ sufficiently large. Here $\scrH_D^1(C,\FG)_\alpha:=\scrH_D^1(C,\FG)\times_{\scrH^1(C,\FG)} \CU_\alpha$, see Remark \ref{Rem_StackGBun},
and $d'=\dim \scrH_D^1(C,\FG)$.

\end{proposition}
\begin{proof}

The first statement follows from Proposition \ref{PropHeckeLocMod}, \cite[Proposition 3.12]{AH_Global} and basic properties of the functor $M^c(-)$.\\
We can take the divisor $D$ sufficiently large such that $\scrH_D^1(C,\FG)_\alpha$ becomes representable by a quasi-projective scheme, e.g. see \cite[Remark 2.9]{AH_Global}. Composing the canonical morphism $M(\scrH_D^1(C,\FG)_\alpha)\to M^c(\scrH_D^1(C,\FG)_\alpha)$ and the morphism $M^c(\scrH_D^1(C,\FG)_\alpha)(d)[2d]\to M^c(\scrH_D^1(C,\FG)_\alpha)$ induced by $\ul\CH_D^\CZ$ gives a morphism in 
$$
\Hom_{DM(k,\BQ)}(M(\scrH_D^1(C,\FG)_\alpha)(d)[2d], M^c(\scrH_D^1(C,\FG)_\alpha\times (C\setminus D)^n)).
$$ 
Note that $\scrH_D^1(C,\FG)_\alpha\to\scrH^1(C,\FG)_\alpha$ is a $\FG_D$-torsor, and since $\FG$ is smooth, by Remark \ref{Rem_StackGBun} and \cite{Kel} $(M(\scrH_D^1(C,\FG)_\alpha)$ is dualizable, therefore the above is isomorphic to $\Hom(\BQ(d+d')[2\cdot(d+d')], M^c(\scrH_D^1(C,\FG)_\alpha\times\scrH_D^1(C,\FG)_\alpha\times C^n))$. Thus we obtain a cycle $\phi_D^\CZ$ in $Ch_{d+d'}(\scrH_D^1(C,\FG)_\alpha\times\scrH_D^1(C,\FG)_\alpha\times (C\setminus D)^n)$, e.g. see \cite[Theorem 8.4]{CD}.

\end{proof}

\begin{remark}
In the above proposition, it is not necessary to take the coefficients in $\BQ$. Namely, when $char~k=0$ (resp. $char~k=p$) one can simply work with $\BZ$ (resp. $\BZ[1/p]$) as the corresponding ring coefficients.    
\end{remark}

\begin{remark}
When $\FG$ is constant for a split reductive group $G$ over $k$, imposing semi-stability condition, in the sense of \ref{Not_Rem1}\ref{Rem-Not e} we get a proper scheme $\scrH^1(C,\FG)^{ss}$, then according to Proposition \ref{PropHeckOverBunG} we can observe that there is a morphism 
$$
\CC_d^{\CL_n^+\FG}(GR_n)\to Ch_{d+d'}(\scrH^1(C,\FG)^{ss}\times_k \scrH^1(C,\FG)^{ss}\times_k C^n)
$$
Here $\CC_d^{\CL_n^+\FG}(GR_n)$ denote the $\BZ$-module generated by closed $d$-equidimensional $\CL_n^+\FG$-equivariant subschemes of $GR_n$.   
\end{remark}

\begin{remark}
When $\FG$ is a constant for a split reductive group $G$ over $k$, there is a conjectural description of the motive with compact support $M^c(\scrH^1(C,\FG))$ of the moduli of $\FG$-bundles according to \cite{BeDh}. For $\FG=\GL_n$, the formula  has been proved inside Voevodsky's motivic categories in \cite{Vicky-Pepin}.
\end{remark}


\begin{definition-remark}[Functoriality]\label{Def-Rem_Funtoriality}
\begin{enumerate}

\item
Consider a morphism $\varrho: \FG' \to \FG$ of algebraic groups. This induces a 1-morphism $\varrho_\ast: \scrH^1(C,\FG')\to\scrH^1(C,\FG)$  and consequently
\begin{equation}\label{Eq_HeckeFunctorial}
\CH_n(C,\FG'):=Hecke_n(C/S,\FG')\to\CH_n(C,\FG):=Hecke_n(C/S,\FG),
\end{equation}
of ind-algebraic stacks. When we further assume that $\FG'\to\FG$ satisfies a condition similar to \ref{EqCondQuotient} (i.e. there is a scheme $Y$ affine and of finite type over $X$ with an action $\FG\times_X Y\to Y$ of $\FG$ and a $\FG$-equivariant open immersion $\FG/\FG'\into Y$), it is quasi-affine and of finite type. Let $\scrE:=End(\scrH^1(C,\FG))$ and $\scrE':=End(\scrH^1(C,\FG'))$. Given morphisms $\varsigma\in\scrE$ and $\varsigma'\in\scrE'$, with $\varsigma\circ\varrho_\ast \cong \varrho_\ast \circ \varsigma$\forget{
which fit into the following diagram
$$
\CD
\scrH^1(C,\FG')@>>>\scrH^1(C,\FG)\\
@V{\varsigma'}VV @V{\varsigma}VV\\
\scrH^1(C,\FG')@>>>\scrH^1(C,\FG),
\endCD
$$
}
then the above morphism \ref{Eq_HeckeFunctorial} induces the following morphism of ind-algebraic stacks
$$
\Sigma_{\varsigma'}(C,\FG') \to \Sigma_\varsigma(C,\FG).
$$
In particular for $\varsigma=\ul\Theta$ and $\varsigma'=\ul\Theta'$ we get a morphism 
$$
\Sigma_{\ul\Theta'}\to \Sigma_\ul\Theta
$$
of  global affine Grassmannians. 
\item
A \emph{$\Sigma\scrH$-datum} is a tuple $(\FG,\CZ,\varsigma)$ consisting of a group scheme $\FG$ over $C$ together with a boundedness condition $\CZ$ and an endomorphism $\varsigma$ in $\scrE(S)$. A morphism $(\FG',\CZ',\varsigma')\to (\FG,\CZ,\varsigma)$ between $\Sigma\scrH$-data is a morphism $\varrho: \FG' \to \FG$ such that $\CZ' \hookrightarrow \Sigma_{\ul\Theta'}\to \Sigma_\ul\Theta$ factors through $\CZ$, and the following commutative diagram 
$$
\CD
\scrH^1(C,\FG')@>\varrho_\ast>>\scrH^1(C,\FG)\\
@V{\varsigma'}VV @V{\varsigma}VV\\
\scrH^1(C,\FG')@>\varrho_\ast>>\scrH^1(C,\FG).
\endCD
$$
  
\forget{
The above morphism can be used to pull boundedness condition $\CZ\subset \Sigma_\ul\Theta$ back and obtain the bound $\varrho^{-1}(\CZ):=\Sigma_{\Theta'}\times_{\Sigma_{\ul\Theta}}\CZ\subseteq \Sigma_{\ul\Theta'}$. }To a morphism $\varrho: (\FG',\CZ',\varsigma')\to (\FG,\CZ,\varsigma)$ as above, one can assign the following 1-morphism of algebraic stacks

$$
\varrho_\ast: \Sigma_{\varsigma'}^{\CZ'}(C,\FG')\to\Sigma_{\varsigma}^\CZ(C,\FG).
$$

\end{enumerate} 
 
\end{definition-remark}  
  
\noindent 
The following proposition describes the local geometry  of the fibers of the family $\Sigma$ over $\scrE$.

\begin{proposition}\label{Prop_LocalModel}

Assume that the group $\FG$ is smooth over $C$. Then there is the following roof of morphisms

\[ 
\xygraph{
!{<0cm,0cm>;<1cm,0cm>:<0cm,1cm>::}
!{(0,0) }*+{\wt\Sigma}="a"
!{(-1.5,-1.5) }*+{\Sigma}="b"
!{(1.5,-1.5) }*+{\Sigma_\Theta \times \scrE,}="c"
"a":_f"b" "a":^{\pi^{loc}}"c"
}  
\]

where

\begin{enumerate}
\item $\wt \Sigma$ is a $\CL_n^+\FG$-torsor over $\Sigma$ under $f$, 
\item
$\pi^{loc}$ is formally smooth.
\end{enumerate}

In addition, for a given $\Sigma\scrH$-datum $(\FG,\CZ,\varsigma)$ the above roof induces a roof of morphisms 

\[ 
\xygraph{
!{<0cm,0cm>;<1cm,0cm>:<0cm,1cm>::}
!{(0,0) }*+{\wt\Sigma_\varsigma^\CZ(C,\FG)}="a"
!{(-1.5,-1.5) }*+{\Sigma_\varsigma^\CZ(C,\FG)}="b"
!{(1.5,-1.5) }*+{\CZ,}="c"
"a":_f"b" "a":^{\pi^{loc}}"c"
}  
\]
in a functorial way.

\end{proposition}

\begin{proof}

Let $\wt \CH$ denote the stack whose category of $T$-valued points parametrizes the tuples $(\ul\CG,\epsilon)$, where $\ul\CG:=(\CG,\CG',\ul s,\tau:\CG\to\CG')$ in $\CH(T)$ and $\epsilon$ is a trivialization of the restriction $\wh\CG'$ of $\CG'$ to the formal neighborhood of $\Gamma_\ul s$. Sending $(\ul\CG,\epsilon)$ to $(\CG,\ul s, \epsilon \circ \tau^{-1})$ gives a map $\wt\CH \to \Sigma_\ul\Theta$. Let $\wt\Sigma$ (resp. $f$) denote the stack (resp. the morphism) defined by the following diagram

$$
\CD
\wt\Sigma @>>>\wt\CH@>>>\Sigma_\ul\Theta\\
@VV{f~~~~~~~\square}V @VVV\\
\Sigma@>>>\CH\\
@VV{~~~~~~~~\square}V @VVV\\
\scrE\times\CY@>>>\CY\times\CY.\\
\endCD
$$

This gives the desired roof of morphisms

\[ 
\xygraph{
!{<0cm,0cm>;<1cm,0cm>:<0cm,1cm>::}
!{(0,0) }*+{\wt\Sigma}="a"
!{(-1.5,-1.5) }*+{\Sigma}="b"
!{(1.5,-1.5) }*+{\Sigma_\Theta \times \scrE,}="c"
"a":_f"b" "a":^{\pi^{loc}}"c"
}  
\]

The morphism $f$ is simply given by forgetting the trivialization $\epsilon$, hence it is a torsor under the group $\CL_n^+\FG$.\\
It remains to justify that $\pi^{loc}:\wt\Sigma \to \Sigma_\Theta \times \scrE$ is formally smooth. Consider the following commutative diagram

\[
\xymatrix {
\ol T\ar[rr]^{(\ol\varsigma,\ol\CG,\ul{\ol s},\ol\tau: \ol\varsigma\ol\CG\to\ol\CG,\ol\epsilon)}\ar[d]_i & & \wt\Sigma\ar[d]^{\pi^{loc}}  \\
T\ar[rr]_{(\CG, \ul s,\phi, \varsigma)}\ar@{.>}[urr]^{\alpha}& & \Sigma_\ul\Theta\times \scrE \;.
}
\]

Where $i: \ol T\to T$ is defined by a nilpotent sheaf of ideals. We need to show that the map $\alpha$ which fits in the above commutative diagram exists. This  question reduces to the fact that a trivialization of a $\FG$-bundle $\wh{\ol{\ul\CG}}$ over infinitesimal neighborhood of $\Gamma_{\ol{\ul s}}$ lifts over $\Gamma_\ul s$. This holds due to smoothness of $\FG$ and the infinitesimal criterion for smoothness. Compare also  proof of \cite[Theorem 3.12]{Ara_LMLocSht}. Part $b)$ follows from part $a)$ and definition of boundedness condition, see Definition-Remark and Definition-Remark \ref{Def_Rem_BC}, and Definition-Remark \ref{Def-Rem_Funtoriality}.

\end{proof}

\subsubsection{Moduli of Higgs bundles}\label{SubsectHiggsBundles}

Let $\FZ$ denote the center of $\FG$. The stack $\scrH^1(C,\FZ)$ operates on $\scrH^1(C,\FG)$ via twisting torsors
$$
\scrH^1(C,\FZ)\times \scrH^1(C,\FG) \to \scrH^1(C,\FG).
$$
It is defined by sending $(\CT,\CG)$ to $\CT\times^{\FZ}\CG$. This induces the following morphism of stacks
\begin{equation}\label{EqH1toEnd}
\scrH^1(C,\FZ)\to \scrE,~\CT\mapsto \varsigma_{\CT}
\end{equation}
(resp. $\scrH^1(C,\FZ)\to \ul{\Hom}:=\ul{\Hom}(\scrH^1(C,\FG)_{(\ul x,\ul P)}^{\ul a-(s)s},\scrH^1(C,\FG))$ obtained by composing \ref{EqH1toEnd} with natural morphism $\scrE\to \ul\Hom$). The pull back of the family $\Sigma^\CZ\to\scrE$ (resp. $\Sigma_{(\ul x, \ul P)}^{\CZ, \ul a-(s)s}\to \ul{\Hom}$ under the above map yields a family $(\Sigma^{\CZ})^\dagger\to\scrH^1(C,\FZ)$

$$
\CD
(\Sigma^{\CZ})^\dagger@>>>\Sigma^\CZ\\
@VVV @VVV\\
\scrH^1(C,\FZ)@>>> \scrE.
\endCD
$$ 

(resp. $\left(\Sigma_{(\ul x, \ul P)}^{\CZ, \ul a-(s)s}\right)^\dagger \to \scrH^1(C,\FZ)$

$$
\CD
\left(\Sigma_{(\ul x, \ul P)}^{\CZ, \ul a-(s)s}\right)^\dagger@>>>\Sigma_{(\ul x, \ul P)}^{\CZ, \ul a-(s)s}\\
@VVV @VVV\\
\scrH^1(C,\FZ)@>>> \ul{\Hom}.)
\endCD
$$

\begin{proposition}\label{PropHiggs^ssisProper}
Let $\FG$ be a constant reductive group, i.e. $\FG=G\times_S C$ for reductive group $G$ over $S$. The morphism $\left(\Sigma_{(\ul x, \ul P)}^{\CZ, \ul a-(s)s}\right)^\dagger\to \scrH^1(C,\FZ)$ is proper, for relevant type $(\ul x,\ul P)$ and stability condition $\ul a$.

\end{proposition}

\begin{proof}
This statement follows from \cite[Proposition~3.12]{AH_Global}, Proposition \ref{PropHeckeLocMod} and Theorem \ref{ThmHeinlothSS}. 

\end{proof}

Let $S=\Spec R$, for a complete dvr $(R,\Fm_R)$, with special point $s$ and generic point $\eta$. Let $\kappa(s):=R/\Fm_R$ denote the residue field and set $\ol S:=\Spec \kappa(s)$. 
\noindent
Let $\CZ$ be a boundedness condition and let $\ol\CZ:=\CZ\times_S\ol S$. 
\noindent
We set $\CY':=\scrH^1(C,\FZ)$. Let $\ol\CY':=\CY'\times_S\ol S$, $\ol\CH:=\CH\times_R \ol S$ and set $\ol{pr}^\leftarrow:=pr^\leftarrow\times_S\ol S$ (resp. $\ol{pr}^\rightarrow:=pr^\rightarrow\times_S\ol S$).  Let $\ol\Sigma$ denote the corresponding family over $\ol\scrE:=\End(\ol\CY)$, see Definition \ref{Def_HR}. We similarly use the notation $\ol\Sigma^\CZ$ for the corresponding bounded family and $\ol\Sigma_{(\ul x, \ul P)}^{\CZ, \ul a-(s)s}$ for the corresponding family with stability conditions of the given type $(\ul x, \ul P)$.

\begin{corollary}\label{CorExtHggs}
Keep the above notation. Let $\CT_0$ be a $\FZ$-bundle in $\ol\CY'(\ol S)$. Twisting by $\CT_0$ defines a morphism $\varsigma_{\CT_0}$ in $\ol\scrE$, and consequently in $\ul\Hom$. Consider the family $\ol\Sigma^{\CZ}\to \ol\scrE$ (resp. $\ol{\Sigma}_{(\ul x, \ul P)}^{\CZ, \ul a-(s)s}\to\ul\Hom$) and let $\ol\CX:=\ol\Sigma_{\varsigma_{\CT_0}}^{\CZ}$ (resp.  $\ol\CX':=\left(\ol{\Sigma}_{(\ul x, \ul P)}^{\CZ, \ul a-(s)s}\right)_{\varsigma_{\CT_0}}$) denote the fiber above $\varsigma_{\CT_0}$. Then\\

a) there is a deformation $\wh\CX$ (resp. $\wh\CX'$) of $\ol\CX$ (resp. $\ol\CX'$) over $S$, \\

b) Let $\ul\nu$ be a $n$-tuple of distinct sections $\nu_i: S\to C$. Assume that $\CZ_{\nu_i}:=\CZ\times_{C,\nu_i} S$ is smooth \forget{(e.g. it is given by a Schubert variety inside affine Grassmannian which is associated with minuscule coweight $\mu_i$)} for every $i$, and let $\wh\CX_\ul\nu'$ denote the corresponding fiber above $\ul\nu:S\to C^n$ of the projection map $\wh\CX'\to C^n$. Then

$$
\Koh^m((\wh\CX_\ul\nu')_\ol\eta,\BQ_\ell)\cong \Koh^m((\wh\CX_\ul\nu')_\ol s,\BQ_\ell).
$$

\end{corollary}

\begin{proof}
The $\FZ$-bundle $\CT_0$ defines a point $\ol S \to \scrH^1(C,\FZ)\to\ul\Hom$. The statement follows from Proposition \ref{PropHiggs^ssisProper}, Proposition \ref{Prop_LocalModel} and Lemma \ref{LemLiftOfCGExist}; see \cite[Proposition~6.0.18]{Wang} for smoothness of $\scrH^1(C,\FZ)$. Note that by Proposition \ref{Prop_LocalModel} we observe that $R\Psi(\ol\BQ_\ell)$ is constant.

\end{proof}

Bellow we address the formally properness of the induced families over $Pic(C/S)$. For the  notion of formally proper morphism see \cite{HLP}.
 
\begin{corollary}
Suppose $C \to S$ admits a section. For a cocharacter $\lambda$ of $\FZ$, the family $(\Sigma^\CZ)^\dagger$ yields a family $(\Sigma^\CZ)_\lambda^\dagger$ over $\Pic(C/S)$. After restricting to $\ol \CU_\alpha$ this gives a formally proper family $(\Sigma_\ol\alpha^\CZ)_{\lambda}^{\dagger}\to \Pic(C/S)$.
\end{corollary}

\begin{proof}
The cocharacter $\lambda$ induces a closed immersion $$\scrH^1(C,\BG_m)\to \scrH^1(C,\FZ).$$ Note that the compositum $\ul\lambda$ of this morphism followed by \ref{EqH1toEnd} is given by sending a $\BG_m$-torsor $\CL$ to the morphism $\ul\CL$, defined by sending $\CG$ to it's twist $\CG\times_{\lambda}^{\BG_m}\CL$ by $\CL$. We may restrict the family $(\Sigma^\CZ)^\dagger\to\scrH^1(C,\FZ)$ to obtain a family $(\Sigma_\ol\alpha^\CZ)_{\lambda}^{\dagger}$ on $\scrH^1(C,\BG_m)$. Since $C\to S$ admits a section, there is an isomorphism $\scrH^1(C,\BG_m)\cong \Pic (C/S) \times B\BG_m$. This gives a morphism $(\Sigma_\ol\alpha^\CZ)_{\lambda}^{\dagger}\to \Pic(C/S)$. The formal properness follows from proposition \ref{PropHiggs^ssisProper} and \cite[Example 4.3.1]{HLP}.
\end{proof}

\begin{remark}[Hitchin morphism for $\FG=\GL_r$]
Let $\FG:=\GL_r$ and let $\lambda$ be the embedding $\BG_m \cong \FZ\to \FG$ of the center. Fix an isomorphism $\CV(-):\scrH^1(C,\FG)\cong Vect_C^r$, where $Vect_C^r$ denotes the stack of vector bundles of rank $r$ over $C$. One can construct a morphism from $(\Sigma^\CZ)_\lambda^{\dagger}$ to an affine bundle over $\scrH^1(C,\FG)\times C^n$ in the following way. Let us take the following variant of the boundedness conditions. For a point $(\CG,\CG',\ul s, \tau)$ in  $Hecke_n(C,\FG)(T)$ we require that $\tau^{-1} \CV(\CG)\subseteq \CV(\CG') \otimes \CO_{C_T}(D_\ul N(\ul s))$, where $D_\ul N(\ul s)$ denotes the relative divisor $\sum_i N_i\cdot \Gamma_{s_i}$. Let $Hecke_n^\ul N(C,\FG)$ (resp. $\Sigma^{\ul N}$) be the moduli stack obtained by imposing this boundedness condition. Then $T$-points of $(\Sigma^\ul N)_\lambda^{\dagger}$ can be described as tuples $(\CV,\ul s, \tau,\CL)$, where $\tau$ is a morphism of vector bundles with $\tau^{-1}\CV\subseteq \CV\otimes \CL(D_\ul N(\ul s))$, and $\CL$ is a line bundle over $C\times T$. Here $\CL(D_\ul N(\ul s)):=\CL\otimes_C \CO(D_\ul N(\ul s))$\forget{ whose cokernel is supported on $D$}. Consider the affine bundle $\CA(-,\ul N)$ over $Pic(C)$, defined by the following functor of points

$$
T \mapsto \{(\CL, \ul s:=(s_i)_i, \ul t:=(t_i)_i);~\CL \in Pic(C)(T),~ s_i \in C(T),~ t_i\in \Koh^0(C_T,\CL^{\otimes i}(i \cdot D_\ul N(\ul s)) \}
$$

\noindent
over $C^n$.\\

A point $\ul\CV:=(\CV,\ul s, \tau,\CL)\in (\Sigma^\ul N)_\lambda^\dagger$, determines a global section $\Tr \tau^{-1}$ of $\CL(D_\ul N(\ul s))$ via the composition $\CO_C \to \check{\CV} \otimes \CV \to \CL(D_\ul N(\ul s))$, where $\check{\CV} \otimes \CV \to \CL(D_\ul N)$ is induced by $\tau^{-1}$ and the first map takes $1$ to $\id_\CV$. Similarly we define $\Tr(\wedge^i \tau^{-1})$. This yields the following map

\begin{equation}\label{Eq_HitchinTypeMap}
(\Sigma^\ul N)_\lambda^{\dagger} \to \CA(-, \ul N).
\end{equation}

\end{remark}

\begin{remark}[The fiber above $\id$ and periodic dynamics]

Let $\varsigma=\id\in \scrE$ (or equivalently consider the trivial bundle $\CL=\CO_C$), let $\Sigma_{\id}$ be the fiber above the identity. Fix a representation $\rho: \FG\to \GL(\CV_0)$, where $\CV_0$ is a vector bundle over $C$. Assume that the induced morphism $\Sigma_{\id}^\CZ(C,\FG)\hookrightarrow Hecke_n(C,\GL_r)$ factors through $Hecke_n^{\ul N}(C,\GL_r)$ for some $\ul N$. Then composing $\Sigma_{\id}^\CZ(C,\FG)\to \Sigma_{\id}^\CZ(C,\GL_r)$ with \ref{Eq_HitchinTypeMap}, for $\CL=\CO_C$, gives

\begin{equation}\label{Eq_HitchinTypeMapII}
\Sigma_{\id}^{\CZ}(C,\FG) \to \CA(\ul N).
\end{equation}
Here $\CA(\CO_C, \ul N):=\CA(\CO_C, \ul N)$ denote the fiber of $\CA(-, \ul N)\to Pic(C)$ over the trivial bundle $\CO_C$.
Assume that $\varsigma^m=\id$ and consider the map $\alpha^m: \Sigma_\varsigma\to\Sigma_{\varsigma^m}$, defined by sending $(\CG,\varsigma\CG,\ul s, \tau: \varsigma\CG\to \CG)$ to $(\CG,\varsigma^m \CG, \ul s, \tau\circ\dots \varsigma^{m-1}\tau)$. This gives $\alpha: \Sigma_\varsigma\to\Sigma_{\id}$. After imposing the boundedness conditions, and composing with the map \ref{Eq_HitchinTypeMap}, we obtain 

$$
\Sigma_{\varsigma}^{\CZ}(C,\FG) \to \CA(\ul N).
$$

\end{remark}

\noindent
\subsubsection{Moduli Of $\FG$-Shtukas}\label{SectModuliG-Sht}

To construct the moduli stack of global $\FG$-shtukas, one needs Frobenius symmetry on $\CY$. To provide this, we have to pass to the formal completions at a fixed prime $p\in\BZ$.

Let $S=\Spec R$ (resp. $\wh S=\Spf R$) for a complete dvr $(R,\Fm_R)$ with special point $s$ and generic point $\eta$. Let $\kappa(s)$ denote the residue field $R/\Fm_R$ and set $\ol S:=\Spec \kappa(s)$. 
\noindent
Let $\FG$ be a smooth affine group scheme over $C$, let $D\subseteq C$ be an effective relative divisor, and let $\CZ$ be a boundedness condition. Set $\ol C:=C\times_S \ol S$, $\ol\FG:= \FG\times_S \ol S$, $\ol\CZ:=\CZ\times_S \ol S$ and $\ol D=D\times_S \ol S$.

In the sequel we explain the analogs for the observations stated in \ref{SubsectHiggsBundles}.

\begin{theorem}\label{Thm_FormalDM} Let $\CW$ be a smooth stack and let $f:\CW\to \CY$ be a morphism of algebraic stacks over $S$. Assume that $\ol\sigma_\CW:=\ol\sigma\times_\CY \ol\CW$ lifts to a morphism $\wh\sigma_\CW:\wh\CW:=\CW\hat{\times}\Spf R \to\wh\CW$. We have the following statements
\begin{enumerate}
\item\label{FA}
There is a formal algebraic stack $\wh\Sigma_{D, \hat{\ul\sigma}}^\CZ$ over $\Spf R$ whose special fiber $\ol\Sigma_{D, \hat{\ul\sigma}}^\CZ$ is Deligne-Mumford and coincides $\nabla_n^{\ol\CZ} \scrH_\ol D^1(\ol C,\ol \FG)\times_\CY \ol\CW$. Here $\nabla_n^{\ol\CZ} \scrH_\ol D^1(\ol C,\ol \FG)$ denotes the moduli of global $\ol\FG$-shtukas.  
\item\label{FA1}

There is a natural morphism $\wh\Sigma_{D', \hat{\ul\sigma}}^\CZ\to \wh\Sigma_{D,\hat{\ul\sigma}}^\CZ$, for $D\subseteq D'$, which is an \'etale morphism of formal algebraic stacks.
 
\item\label{FA'}
Assume that $\CW$ is proper. Then the restriction $\wh\Sigma_{D, {\hat{\ul\sigma}}, \ol\alpha}^\CZ$ of $\wh\Sigma_{D, {\hat{\ul\sigma}}}^\CZ$ to $\ol \CU_\alpha \times_\CY \wh\CW$ is algebraizable for large enough $D$, i.e. there is $\Sigma_{D, {\hat{\ul\sigma}}, \ol\alpha}^\CZ$ over $\Spec R$ with $\wh\Sigma_{D, {\hat{\ul\sigma}}, \ol\alpha}^\CZ=\Sigma_{D, {\hat{\ul\sigma}, \ol\alpha}}^\CZ\wh\times_R \Spf R$.

\end{enumerate}

\end{theorem}

\begin{proof}

\ref{FA} $\CY_D\to\CY$ is a $\FG_D$-torsor, where $\FG_D$ denote the Weil restriction of scalars $R_{D/S}\FG$. Let $\wh\CY:=\CY\times_S \wh S$ and $\wh\CY_D:=\CY_D\times_S \wh S$. Let us set $\wh\CW_D:=\wh\CW\hat{\times}\wh\CY_D$. The map $\wh\sigma_\CW$ lifts to a morphism $\hat{\ul\sigma}:=\hat{\ul\sigma}_D=\wh\CW_D\to\wh\CW_D$. Consider $\wh\CH_D^\CZ:= \CH_D^{\CZ} \wh\times \wh S$ as a family over $\wh\CY_D\wh\times\wh\CY_D$ and let $\hat{H}:=\wh\CH_D^\CZ\times_{\wh\CY_D\times\wh\CY_D}(\wh\CW_D\times\wh\CW_D)$. Define

$$
\wh\Sigma_{D, \hat{\ul\sigma}}^\CZ:=\Sigma(\ul{\hat{H}})_{\ul{\wh\sigma}}.
$$ 

Here $\ul{\hat{H}}$ denote the $HR$-tuple $(\hat{H}, char: \hat{H}\to \wh{C}^n:=C^n\wh\times \wh{S}, pr^{\leftarrow}, pr^{\rightarrow}: \hat{H} \rightrightarrows\wh\CW_D)$. By construction the special fiber $\ol\Sigma_{D, \hat{\ul\sigma}}^\CZ$ coincides $\nabla_n^{\ol\CZ} \scrH_\ol D^1(C,\ol \FG)\times_\CY \ol\CW$; see \cite{AH_Global} and \cite{AH_LM}. In particular it is Deligne-Mumford; see \cite[Theorem~3.1.7]{AH_LM}.

\ref{FA1} Flatness of this morphism can be checked over the special fiber according to \cite[Lemme~11.3.10.1]{EGA}; see also Definition-Remark \ref{Def_Rem_BC} b). Now the statement follows from \cite[Theorem~3.15]{AH_Global} and \cite[Lemma 1.2]{BL}.\\
\ref{FA'}
Since $\CZ$ is proper, see Definition-Remark~\ref{Def_Rem_BC} b), by Proposition \ref{PropHeckeLocMod} we observe that $\CH_\ol\alpha^\CZ \to  \CY_\ol\alpha:=\ol\CU_\alpha$, and therefore $\CH_\ol\alpha^\CZ \to  \CY_\ol\alpha \times  \CY_\ol\alpha$, are proper. 
Consequently $\wh\Sigma_{D, \hat{\ul\sigma}, \ol\alpha}^ \CZ$ is proper over $\wh\CW_{\ol\alpha}:=\wh\CW\times_\CY\CY_\ol\alpha$. Notice that $\wh\Sigma_{D, \hat{\ul\sigma}, \ol\alpha}^\CZ$ is equipped with an ample line bundle (i.e. a system of line bundles $\CL_n$ on $(\wh\Sigma_{D, \hat{\ul\sigma}, \ol\alpha}^\CZ)_n:=\wh\Sigma_{D, \hat{\ul\sigma}, \ol\alpha}^\CZ \times R/\Fm_R^n$) which is inherited from an ample line bundle on $\CH_{D,\ol\alpha}^\CZ$ for large enough $D$. The existence of ample line bundle on $\CH_{D,\ol\alpha}^\CZ$ is a consequence of Proposition \ref{PropHeckeLocMod} and \cite[Thorem~5.0.14]{Wang}. Now the statement follows from Grothendieck's algebraization theorem;  \cite[III,~Thm. 5.4.5]{EGA}

\end{proof}

\begin{proposition}\label{ProExExissForAffChars}
 Assume that $\CW\to \CY$ is a smooth affine chart, then there is a lift $\wh\Sigma_{D, \ul\sigma}^\CZ$ of $\nabla_n^{\ol\CZ} \scrH_\ol D^1(\ol C,\ol \FG)\times_\CY \ol\CW$ as in Theorem \ref{Thm_FormalDM} \ref{FA}.
\end{proposition}

\begin{proof}
As $\FG$ is smooth we see that $\CY$, and thus $\CW$, are smooth. Since $\CW$ is in addition affine, we see by infinitesimal criterion for smoothness that the Frobenius $\ol\sigma_\ol\CW$ lifts, and gives $\wh\sigma: \wh\CW\to\wh\CW$.   
\end{proof}

Fix a n-tuple $\ul\nu:=(\nu_i)$ of disjoint characteristic sections on $C$. That is $\ul\nu: \Spf R \to C^n$. Let $\CZ_\ul\nu:=(\CZ_{\nu_i})$ denote the associated tuple of local bounds corresponding to $\CZ$ at the places $\nu_i$, see \cite[Subsection 4.3]{AH_LM}.

\begin{proposition}\label{PropKohSpec&generic}
Keep the notation in theorem \ref{Thm_FormalDM}. Furthermore assume that $\CW$ is proper and $\CZ_{\nu_i}$ is smooth (e.g. it comes from minuscule coweights) for every $i$, then there is an isomorphism

$$
\Koh^q((\wh\Sigma_{D,\hat{\ul\sigma},\ul\nu}^\CZ)_\ol\eta,\BQ_\ell)\tilde{\to} \Koh^q(\nabla_n^{\ol\CZ_\ul\nu} \scrH_\ol D^1(C,\ol \FG)\times_\CY \ol\CW
,\BQ_\ell).
$$
In particular the cohomology of the generic fiber $(\wh\Sigma_{D,\ul\sigma,\ul\nu}^\CZ)_\ol\eta$ is independent of the choice of the lift $\ul\sigma$. 
\end{proposition}

\begin{proof}

We may take $D$ enough large such that $\wh\Sigma_{D,\hat{\ul\sigma}}$ becomes algebraizable. The projection morphism $\CH^\CZ\to \CY$ is proper, see proof of theorem \ref{Thm_FormalDM} c). Since $\CW$ is proper, we observe that $\CH^\CZ\times_{\CY\times\CY}(\CW\times\CW) \to \CW\times\CW$ is proper and hence $\wh\Sigma_{D,\hat{\ul\sigma}}$ is proper. Now the isomorphism follows from Lemma \ref{LemLiftOfCGExist} and Proposition \ref{Prop_LocalModel} which implies that $R\Psi(\ol\BQ_\ell)$ is constant; see also Lemma \ref{LemLocalModelI} below.

\end{proof}

\begin{remark}
For the case $\ell=p$, and for the tuple of smooth bounds $\CZ_\ul\nu$, it can be seen that the cohomology $\Koh^q((\wh\CX_{\hat{\ul\sigma}})_\eta,\BQ_\ell)$ of $\wh\CX_\ul\sigma=(\wh{\Sigma}_{D,\hat{\ul\sigma},\ul\nu})_\eta$ is also independent of the choice of the lift $\hat{\ul\sigma}$. This follows from the spectral sequence  

$$
\Koh^p(\FX_\ol s, i^\ast R^q j_\ast \BQ_p)\Rightarrow \Koh^{p+q}(\FX_\ol\eta,\BQ_p),
$$

\noindent
see for example \cite{BK}, corresponding to the following diagram of formal stacks

$$
\CD
\FX_\ol\eta @>>> \FX @<<< \FX_\ol s\\
@VVV @VVV @VVV\\
\ol\eta @>>>\ol S @<<< \ol s.
\endCD
$$
\noindent
The bar in the above notation indicates corresponding algebraic or integral closures.

\end{remark}

\begin{question}

\begin{enumerate}

Fix a global boundedness condition $\CZ$. Consider the stack $\CY_{(\ul x, \ul P)}^{\ul a-ss}:=\scrH^1(C,\FG,\ul x, \ul P)^{\ul a-ss}$ of $\ul a$-semistable flagged principal $\FG$-bundles of type $(\ul x, \ul P)$. 

It is natural to ask 

\item[-]
whether one can take $\CW$ to be $\CY_{(\ul x, \ul P)}^{\ul a-ss}$ (i.e. if the Frobenius lifts over $\wh\CW=\CY_{(\ul x, \ul P)}^{\ul a-ss}\wh\times \Spf R$)? and second,

\item[-]
does the cohomology remain independent of the choice of the lift $\hat{\ul\sigma}$ in the global situation, namely for the stack $\wh\Sigma_{\hat{\ul\sigma}}^\CZ$?
Note that according to the lemma \ref{LemLocalModelI}, the moduli stack $\wh\Sigma_{\hat{\ul\sigma}}^\CZ$ is flat over $C^n$, provided that $\CZ$ is flat over $C^n$. Furthermore, as $\pi:\CX_{\hat{\ul\sigma}}:=(\wh\Sigma_{\ul\sigma}^{\CZ_\ul\nu})_{(\ul x, \ul P)}^{\ul a-s} \to C^n$ is proper, the higher direct image $R^i~{\pi}_\ast \CF$ is constructible and we have the following convergence of Leray spectral sequence

$$
\Koh^i(C^n,R^i\pi_{\ast}\CF)\Rightarrow\Koh^{i+j}(\CX_{\hat{\ul\sigma}}, \CF).
$$
Moreover we have
$$
(R^i~{\pi}_\ast \CF)_y\cong \Koh^i((\CX_{\hat{\ul\sigma}})_y, \CF_{|y})
$$
e.g. see \cite[Ch. VI Cor. 2.7]{Milne}. 
\end{enumerate}

\end{question}

\begin{lemma}\label{LemLocalModelI}
Keep the notations in theorem \ref{Thm_FormalDM} and assume that $\CW\to\CY$ is \'etale (resp. smooth). Let $\CZ$ be a boundedness condition in the sense of \ref{Def_Rem_BC} b). For a point $y$ in $\wh{\Sigma}_{D,\hat{\ul\sigma}}$ there exist an \'etale neighborhood $U_y$ of $y$ and a roof
\[ 
\xygraph{
!{<0cm,0cm>;<1cm,0cm>:<0cm,1cm>::}
!{(0,0) }*+{U_y}="a"
!{(-1.5,-1.5) }*+{\wh\Sigma_{D, \hat{\ul\sigma}}^\CZ}="b"
!{(1.5,-1.5) }*+{\wh\CZ,}="c"
"a":_{\text{\'et}}"b" "a":^{\text{\'et (resp. smooth)}}"c"
}  
\]
of \'etale morphisms. Here $\wh\CZ:=\CZ\times_R \Spf R$. 
\end{lemma}

\begin{proof}

By theorem \ref{Thm_FormalDM} b) we can forget $D$-level structure. Let $y'$ be the image of $y$ in $C^n\times\CY$ under the projection sending $(\CG,\CG',\ul s,\tauGlob)$ to $(\ul s,\CG')$. According to Proposition~\ref{PropHeckeLocMod}, we may take an \'etale neighborhood $U\to C^n\times_S \CY$ of $y'$, such that the restriction $U'$ of $\CH$ to $U$ and the restriction $U''$ of $\Sigma_\ul\Theta \times \CY$ to $U$ become isomorphic. Now, set $\wt U:=U' \times_{\CH}\wh\Sigma_\sigma^\CZ$. We deduce the following roof of morphisms

\[ 
\xygraph{
!{<0cm,0cm>;<1cm,0cm>:<0cm,1cm>::}
!{(0,0) }*+{\wt U}="a"
!{(-1.5,-1.5) }*+{\wh\Sigma_{\hat{\ul\sigma}}^\CZ}="b"
!{(1.5,-1.5) }*+{\wh\CZ}="c"
"a":_{\text{\'et}}"b" "a":^{\phi}"c"
}  
\]

It remains to check that $\phi$ is an \'etale morphism (resp. smooth). To see this, consider the morphism $\ol\phi:\ol U := U\times_R R/\Fm_RR \to \ol \CZ:=\wh\CZ \times_R R/\Fm_R R$.  The morphism $\ol\phi$ is \'etale according to \cite[Theorem~3.2.1]{AH_LM}. Now as in the proof of \ref{Thm_FormalDM} \ref{FA1} we may argue by \cite[IV, Lemme~11.3.10.1]{EGA} and \cite[Lemma 1.2]{BL}.

\end{proof}

\section{The Hecke stack over the moduli of G-Shtukas}\label{Sect The Hecke stack over the moduli of G-Shtukas}

In this subsection we discuss another sample of the construction we described in section \ref{SectGeneralConstruction}. Namely, we consider the Hecke stack over the moduli stack $\nabla_n\scrH^1(C,\FG)$ (see theorem \ref{Thm_FormalDM} for the notation) of $\FG$-shtukas. Throughout this subsection we let $S:=\Spec\BF_q$.

\begin{definition}\label{DefHeckeoverNablaH}

Fix integers $m$, $n$ and $n'$. Set $\nabla_n\scrH:=\nabla_n\scrH^1(C,\FG)$ and $\nabla_{n'}\scrH:=\nabla_{n'}\scrH^1(C,\FG)$.

\begin{enumerate}

\item
Define the algebraic stack $\wt\CH_{n,m,n'}$ as the stack whose points over a scheme $T$ over $S$ consists of tuples $(\ul\CG,\ul\CG',\ul c,\phi)$, consisting of the following data

\begin{enumerate}

\item 

$\ul\CG:=(\CG,\ul s,\tau)$ in $\nabla_n\scrH(T)$ and $\ul\CG':=(\CG',\ul s',\tau')$ in $\nabla_{n'}\scrH(T)$,

\item

$m$-tuple $\ul c:=(c_i)$ in $C^m(T)$,

\item

a commutative diagram 

$$
\CD
\sigma^\ast\CG'@>{\tau'}>>\CG'\\
@V{\sigma^\ast\phi}VV @VV{\phi}V\\
\sigma^\ast\CG@>{\tau}>>\CG,\\
\endCD
$$ 

where $\phi$ is defined over $C_T\setminus \Gamma_\ul c$ and the diagram is defined after restricting to $C_T\setminus\Gamma_\ul s \cup \Gamma_{\ul s'} \cup \Gamma_\ul c$. 
\end{enumerate}

\item
The moduli stack $\wt\CH_{n,m,n'}$ is fibered over $\nabla_n\scrH$ (resp. $\nabla_{n'}\scrH$) through projections to the first (resp. second) factor. We have the following map

\begin{equation}\label{Eq_HeckeNablaHmorphtoI}
(pr^\leftarrow,pr^\rightarrow, char \circ pr^\leftarrow, char \circ pr^\rightarrow, \wt{char}): \wt\CH_{n,m,n'}\to \nabla_n\scrH\times \nabla_{n'}\scrH\times C^n \times C^{n'} \times C^m,
\end{equation}

given by sending the $T$-point $(\ul\CG:=(\CG,\ul s,\tau),\ul\CG':=(\CG',\ul s',\tau'),\ul c,\phi)$ in $\wt\CH_{n,m,n'}(T)$ to $(\ul\CG,\ul\CG',\ul s, \ul s', \ul c)$.

\item

There is an obvious projection

$$
\CH_{n,m,n'}\to Hecke_m(C,\FG),
$$

defined by sending $(\ul\CG,\ul\CG',\ul c,\phi)$ to $(\CG,\CG',\ul c,\phi)$. Set 
$$
\CH_{n,m,n'}^{\CZ,\wt\CZ,\CZ'}:=\CH_{n,m,n'}\times_{\nabla_n\scrH\times \nabla_{n'}\scrH\times Hecke_m(C,\FG)}\left(\nabla_n^\CZ\scrH\times \nabla_{n'}^{\CZ'}\scrH\times Hecke_m^{\wt\CZ}(C,\FG) \right),
$$
Here $\CZ$ (resp. $\wt\CZ$, resp. $\CZ'$) is a bound in $GR_n$ (resp. $GR_m$, resp. $GR_{n'}$).

\item

Let $\Hom_{n,n'}$ denote the stack $\Hom(\nabla_n\scrH,\nabla_{n'}\scrH)$  and let $$\ul\CH_{n,m,n'}:=(\CH_{n,m,n'},\wt{char}, pr^\leftarrow, pr^\rightarrow).$$ 
Let $\Sigma_{n,m,n'/\nabla\scrH}:=\Sigma(\ul\CH_{n,m,n'})$ denote the corresponding family over $\Hom_{n,n'}$. For $n=n'$, set $\scrE_n=\Hom_{n,n'}$.
\item
Let $D\subseteq C$ be a closed subscheme. One can equip the Hecke stack $\CH_{n,m,n'}$ (resp. $\CH_{n,m,n'}^{\CZ,\wt\CZ,\CZ'}$) with a $D$-level structure in an obvious sense. We denote the resulting moduli stack by $\CH_{D,n,m,n'}$ (resp. $\CH_{D,n,m,n'}^{\CZ,\wt\CZ,\CZ'}$).

\end{enumerate}

\end{definition}

\begin{remark}
Let $I$ (resp. $J$, resp. $K$) be an index set with cardinality $n$ (resp. $n'$, resp. $n+n'$) and fix a bijection $\iota: K\to I \sqcup J$. To this, one assigns the following closed immersion of stacks

$$
C^{n'}\times\nabla_n\scrH^1(C,\FG)\to \nabla_{n+n'}\scrH^1(C,\FG),
$$ 

which is defined by sending $(s_j)_{j\in J}\times (\CG, (s_i)_{i\in I},\tau)$ to $(\CG, (s_k)_{k\in K},\tau)$, where $s_k:=s_{\iota k}$.

\end{remark}

\begin{proposition}
Fix a $n$-tuple $\ul\nu:=(\nu_i)$ of closed points $\nu_i$ in $C$. Let $A_\ul\nu$ be the completion of the stalk of $\CO_{C^n}$ at the point $\ul\nu$. Set $\CH_m:=\CH_{n,m,n}$, and let $\CH_{m,\ul\nu}:=\CH_m\times_{C^n} \Spf A_\ul\nu$ (resp. $\nabla_n\scrH_\ul\nu:=\nabla_n\scrH\times\Spf A_\ul\nu$). The projection map $\CH_m \to \nabla_n\scrH_\ul\nu\times \nabla_n\scrH_\ul\nu \times C^m$
is formally unramified.
\end{proposition}

\begin{proof}
This follows from rigidity of quasi-isogenies between $\FG$-shtukas, see \cite[Proposition 5.9]{AH_Local}. 
\end{proof}

\begin{remark}
For an ind-scheme $X^\centerdot=\dirlim[i\in I] X^i$, one defines the Chow group $Ch_d(X^\centerdot)$ as the direct limit of $Ch_r(X^i)$ under $\iota_{i,j_\ast}$, where $\iota_{i,j }:X^i\to X^j$ denotes the closed immersion $X^i\to X^j$ for $i<j$.  For construction of the Chow groups of algebraic spaces, flat pull back, proper push forward and further properties see \cite[Chapter 80, Tag 0EDQ]{Stacks}. Concerning the construction of Chow groups for Deligne-Mumford stacks, e.g. see \cite{Vis} and \cite{Kre}, one similarly defines the (homological) Chow group $Ch_d(\CX^{\centerdot})$, for ind-Deligne-Mumford stack $\CX^{\centerdot}=\dirlim[i\in I] \CX^i$.  Recall that there is an ind-Deligne-Mumford structure $\nabla_n\scrH_D=\dirlim[\ul\omega]\nabla_n^\ul\omega\scrH_D^1(C,\FG)$, where ${\ul\omega}$ lies in a countable set of cocharacters of $SL(\CV)$ for some vector bundle $\CV$; see \cite[Theorem 3.15]{AH_Global}.\\

\end{remark}

\begin{theorem}\label{ThmCyclesOnNablaH}

Fix a boundedness condition $\CZ$, and an $\alpha$ as in Remark \ref{Rem_StackGBun}. Let $\nabla_n^\CZ\scrH_{D,\alpha}:=\nabla_n^\CZ\scrH_{D}\times_{\nabla_n^\CZ\scrH} \CU_\alpha$. Let $\CC_\wt d^{\CL_m^+\FG}(GR_m)$ denote the free $\BZ$-module generated by $\CL_m^+\FG$-equivariant closed equidimensional subschemes of global affine Grassmannian $GR_m$, which are flat over $C^m$ and of dimension $\wt d$. An element $\wt\CZ$ in $\CC_\wt d^{\CL_m^+\FG}(GR_m)$ induces a morphism 
$$
M^c(\nabla_n^\CZ\scrH_{D, \alpha})(\wt d)[2(\wt d)] \to M^c(\nabla_{n+m}^{\CZ'}\scrH_D)\otimes M(C)^{\otimes m},
$$ 
for sufficiently large boundedness condition $\CZ'$. Similarly it induces a morphism in 
$$
\Hom_{\text{DM}_{gm}}(M^c(X_{D,n,\alpha}^\CZ)(\wt d)[2\wt d], M^c(X_{D,n+m}^{\CZ'})\otimes M(C\setminus D)^{\otimes m}),
$$ 
here $X_{D, n,\alpha}^\CZ$ (resp. $X_{D, n+m}^{\CZ'}$) denotes the coarse moduli space for $\nabla_n^\CZ\scrH_{D,\alpha}$ (resp. $\nabla_{n+m}^{\CZ'}\scrH_D$).

\end{theorem}

\begin{proof}
As the statement with $D$-level structure follows similarly, we just explain the situation where there is no level structure.
Note first that for any boundedness condition $\CZ$, the moduli stack $\CX_n^\CZ:=\nabla_n^\CZ\scrH$ (resp. $\CX_{n,\alpha}^\CZ:= \nabla_n^\CZ\scrH_{\alpha}$) is Deligne-Mumford. It is separated and its inertia is finite over $\nabla_n^\CZ\scrH$ (resp. $\nabla_n^\CZ\scrH_{\alpha}$); see \cite[theorem 3.15]{AH_Global}, \cite[theorem 3.1.7]{AH_LM} and \cite[Corollary~3.16]{AH_Global}. Therefore by Keel-Mori's theorem, \cite{Conrad}, it admits a coarse moduli space $X_n^\CZ$ (resp. $X_{n,\alpha}^\CZ$).

Let $\wt\CZ$ be an irreducible cycle in $\CC_\wt d^{\CL_m^+\FG}(GR_m)$. Let $\CH_{(n,m,n+m),\alpha}^{-,\wt \CZ,-}$ denote the substack of $\CH_{n,m,n+m}$ defined by restricting $\CH_{n,m,n+m}$ to $\CU_\alpha$ and then imposing the boundedness condition $\wt\CZ$ to the universal isomorphism $\phi_u$ of the universal tuple $(\ul\CG_u,\ul\CG_u',\ul c_u, \phi_u:\CG_u'\tilde{\to}\CG_u)$.\\
The family $\CH_{(n,m,n+m),\alpha}^{-,\wt\CZ,-}\to \nabla\scrH_\alpha$ is locally on the base isomorphic to the family $\CH_{m}^\wt\CZ\to \scrH^1(C,\FG)$. In particular it is flat of relative dimension $\wt d$; see Proposition \ref{PropHeckeLocMod}. Set

$$
\CH_{(n,m,n+m),\alpha}^{\CZ,\wt\CZ,-}:=\CH_{(n,m,n+m),\alpha}^{-,\wt\CZ}\times_{\nabla\scrH}\nabla_n^\CZ\scrH.
$$
By \cite[Theorem 3.15]{AH_Global} the morphism $pr^\rightarrow: \CH_{(n,m,n+m),\alpha}^{\CZ,\wt\CZ,-}\to \nabla_{n+m}\scrH$ factors through $\CX_{n+m}^{\CZ'}:=\nabla_{n+m}^{\CZ'}\scrH$ for some sufficiently large $\CZ'$. Therefore we get the following roof

\[ 
\xygraph{
!{<0cm,0cm>;<1cm,0cm>:<0cm,1cm>::}
!{(0,0) }*+{\CH_{(n,m,n+m),\alpha}^{\CZ,\wt\CZ,-}}="a"
!{(-1.5,-1.5) }*+{\CX_\alpha^{\CZ}}="b"
!{(1.5,-1.5) }*+{\CX_{n+m}^{\CZ'}\times_S C^n.}="c"
"a":"b" "a":"c"
}  
\]
which induces a morphism $M^c(X_\alpha^\CZ)(\wt d)[2 \wt d]=M^c(\CX^{\CZ})(\wt d)[2 \wt d)] \to M^c(\CX^{\CZ'})=M^c(X^{\CZ'})$.\\ 

\end{proof}

\begin{remark}

Note in particular that when $\nabla_n^\CZ\scrH_D$ is proper we have $$\Hom_{\text{DM}_{gm}}(M^c(X_{n,\alpha}^\CZ)(\wt d)[2\wt d], M^c(X_{n+m}^{\CZ'}))= Ch_{d+\wt d}(X^\CZ\times_S X^{\CZ'}).$$ 
This follows from the above discussion and \cite[Chapter 5, Proposition 4.2.3]{VSF}. Note further that by Zariski's main theorem we observe that $\CX^{\CZ}$ is proper and thus $M^c(\CX^{\CZ})=M(\CX^{\CZ})$. For examples of such cases see \cite{Lau}. 

\end{remark}

\begin{corollary}
Keep the notation in \ref{ThmCyclesOnNablaH}. For $D$ sufficiently large, the element $\wt \CZ$ in $\CC_\wt d^{\CL_m^+(GR_m)}$ induces a cycle $\Phi^\wt\CZ$ in 
$$
A_{d,0}(\nabla_n^{\CZ}\scrH_{D,\alpha},\nabla_{n+m}^{\CZ'}\scrH_{D,\alpha} \times_S (C\setminus D)^m).
$$
Furthermore when $\wt\CZ$ comes from minuscule coweights $\ul\mu:=(\mu_i)_i$, it induces an element in $Ch_{d+\wt d}(\nabla_n^{\CZ}\scrH_D^\circ\times_S\nabla_{n+m}^{\CZ'}\scrH_D \times_S (C\setminus D)^m)$. Here $\nabla_n^{\CZ}\scrH_D^\circ:=\nabla_n^{\CZ}\scrH_D\times_{C^n}(C^n\setminus \Delta^{big})$, where  $\Delta^{big}$ denotes the big diagonal in $C^n$, and $d=\dim \CZ$.
\end{corollary}

\begin{proof}

First we can take $D$ sufficiently large that $\nabla_n^{\CZ}\scrH_{D,\alpha}\times_S\nabla_{n+m}^{\CZ'}\scrH_{D,\alpha}$
becomes representable by a quasi-projective variety. According to \ref{ThmCyclesOnNablaH}, $\wt\CZ$ induces a morphism 
$$
M^c(\nabla_n^{\CZ}\scrH_{D,\alpha})(\wt d)[2\wt d]\to M^c(\nabla_{n+m}^{\CZ'}\scrH_{D,\alpha}\times (C\setminus D)^m).
$$
After composing the canonical morphism $M(\nabla_n^{\CZ}\scrH_{D,\alpha}) \to M^c(\nabla_n^{\CZ}\scrH_{D,\alpha})$, e.g. see \cite[Exampl~16.2]{MVW}, with the above morphism we obtain

$$
M(\nabla_n^{\CZ}\scrH_{D,\alpha})(\wt d)[2\wt d]\to M^c(\nabla_{n+m}^{\CZ'}\scrH_{D,\alpha}\times (C\setminus D)^m)
$$

Note further that $$\Hom(M(\nabla_n^{\CZ}\scrH_{D,\alpha})(\wt d)[2 \wt d], M^c(\nabla_{n+m}^{\CZ'}\scrH_{D,\alpha}\times (C\setminus D)^m))\cong A_{\wt d,0}(\nabla_n^{\CZ}\scrH_{D,\alpha} , \nabla_{n+m}^{\CZ'}\scrH_{D,\alpha}\times (C\setminus D)^m);$$ see \cite[Theorem 8.4]{CD}. Finally when $\CZ$ comes from minuscul coweights, the restriction $\CZ\times_{C^n} (C^n\setminus \Delta^{big})$ is smooth, and therefore $\nabla_n^{\CZ}\scrH_{D,\alpha}^\circ$ is smooth and of dimension $d=n+\sum_i \mu_i$; see \cite[Theorem 3.2.1]{AH_LM}. Hence by duality, see \cite[Theorem 16.24]{MVW} and see also \cite[Theorem 5.4.20]{Kel}, we deduce 
$$
A_{d,0}(\nabla_n^{\CZ}\scrH_{D,\alpha}^\circ, \nabla_{n+m}^{\CZ'}\scrH_{D,\alpha}))=A_{d+\wt d,0}(S, \nabla_n^{\CZ}\scrH_{D,\alpha}^\circ\times_S\nabla_{n+m}^{\CZ'}\scrH_{D,\alpha})
$$
$$
~~~~~~~~~~~~~~~~~~~~~~~~~~~~~~~~~~=Ch_{d+\wt d}(\nabla_n^{\CZ}\scrH_{D,\alpha}^\circ\times_S\nabla_{n+m}^{\CZ'}\scrH_{D,\alpha})
$$

\end{proof}

\begin{remark}[global Rapopor-Zink spaces]
Let $\CZ$ be a boundedness condition in $GR_n$. Let $\CH_{m/\nabla\scrH}^{\CZ,-,-}:=\CH_{m/\nabla\scrH}\times_{\nabla\scrH,pr^\leftarrow}\nabla^{\CZ}\scrH$. Let $\Sigma_{m/\nabla\scrH}^{\CZ,-,-}:=\Sigma(\CH_{m/\nabla\scrH}^{\CZ,-,-},pr^\leftarrow,pr^\rightarrow)$ and view it as a family over $\Hom(\nabla^{\CZ}\scrH,\nabla\scrH)$, according to the construction described in Section \ref{SubsectTwoFiberBundles}. Fix a global $\ul\CG_0$ and let $\ul\Theta_{\ul\CG_0}$ denote the corresponding constant morphism in $\scrE$. The fiber $\Sigma_{m/\nabla\scrH,\ul\Theta_{\ul\CG_0}}^{\CZ,-,-}$ above $\ul\Theta_{\ul\CG_0}$ is called \emph{global Rapoport-Zink space} corresponding to $\CZ$ and $\ul\CG_0$. Note that there is a natural projection morphism 

 $$
 \Sigma_{m/\nabla\scrH,\ul\Theta_{\ul\CG_0}}^{\CZ,-,-}\to \Sigma_{\Theta_{\CG_0}}
 $$

defined by sending $(\ul\CG:=(\CG,\tau_\CG),\ul c ,\phi)$ to $(\CG,\ul c, \phi:\CG|_{C_T\setminus\Gamma_{\ul c}}\to \CG_0|_{C_T\setminus\Gamma_{\ul c}})$; see definition \ref{DefBDGr}.
\end{remark}

\begin{remark}

Consider the morphism $\nabla\scrH\to \nabla\scrH$ induced by Frobenius $\sigma:\scrH^1(C,\FG)\to \scrH^1(C,\FG)$. We again denote this morphism by $\sigma$ and we let $\Sigma_{n/\nabla\scrH, \sigma}$ denote the fiber above $\sigma$. There is a morphism

$$
\Psi_0:\nabla\scrH\to\Sigma_{n/\nabla\scrH, \sigma},
$$ 
which is defined by sending $\ul\CG:=(\CG,\ul s,\tau)$ to $(\ul\CG, \sigma^\ast\ul\CG, \ul s, \tau: \sigma^\ast\ul\CG\to\ul\CG)$.

On the other hand there is a morphism

$$
\Psi_\ell: \Sigma_{/\nabla\scrH,\sigma}\to \Sigma_{n/\nabla\scrH,\sigma^\ell},
$$
defined by sending $(\ul\CG,\sigma^\ast\ul\CG,\ul c,\phi)$ to $(\ul\CG,(\sigma^\ell)^\ast\ul\CG,\ul c,(\sigma^{\ell-1})^\ast\phi\circ\dots\circ\sigma^\ast\phi\circ\phi)$. In particular for a finite field extension $k/\BF_q$ of degree $\ell$, composing $\Psi_\ell$ with  $\Psi_0$ gives

$$
\nabla\scrH(k)\to \Sigma_{n/\nabla\scrH,\sigma^\ell}=\Sigma_{n/\nabla\scrH,\id}(k).
$$

The above map sends $\ul\CG$ to the Frobenius isogeny $\Phi_{\ul\CG}\in QEnd(\ul\CG)$.
Note that, when $\FG=\GL_n$, the coefficients of the minimal polynomial of $\Phi_{\ul\CG}\in QEnd(\ul\CG)$ determines the quasi-isogeny class of $\ul\CG$; see \cite[Section 5]{CMot}.
\end{remark}

\Verkuerzung
{
\vfill

\begin{minipage}[t]{0.35\linewidth}
\noindent
Esmail Arasteh Rad\\
Universit\"at M\"unster\\
Mathematisches Institut \\
erad@uni-muenster.de
\\[1mm]
\end{minipage}

}


{\small
\begin{thebibliography}{GHKR2}
\addcontentsline{toc}{section}{References}






\bibitem[Ao]{Ao} M.~ Aoki \emph{Hom stacks}. preprint available as \href{http://arxiv.org/pdf/math/0503358.pdf}{arXiv:0503358}


\bibitem[Ara]{Ara_LMLocSht} E. Arasteh Rad, \emph{ Rapoport-Zink Spaces For Local $\mathbb{P}$-Shtukas and Their Local Models}, preprint 2018. 32 pages, Available at https://arxiv.org/pdf/1807.03301.pdf

\bibitem[AraHab]{AH_LM} E.~Arasteh Rad and S.~Habibi. \emph{Local Models For the Moduli Stacks of Global G-Shtukas}, Mathematical Research Letters, Vol. 26, No. 2 (2019), pp. 323-364, preprint available as \href{http://arxiv.org/abs/1605.01588v3}{arXiv:1605.01588v3}.


\bibitem[AraHar1]{AH_Local} E.~Arasteh Rad, U.~Hartl: \emph{Local $\BP$-shtukas and their relation to global $\FG$-shtukas}, Muenster J.~Math (2014); also available as \href{http://arxiv.org/abs/1302.6143}{arXiv:1302.6143}.

\bibitem[AraHar2]{AH_Global} E.~Arasteh Rad, U.~Hartl: \emph{Uniformizing the moduli stacks of global $\FG$-Shtukas},  Int. Math.
Res. Not. (2019), https://doi.org/10.1093/imrn/rnz223, in press, also available as \href{http://arxiv.org/abs/1302.6351}{arXiv:1302.6351}.


\bibitem[AraHar3]{CMot} E.~Arasteh Rad and U. Hartl \emph{Category of $C$-Motives over Finite Fields}, Journal of Number Theory 2020 https://doi.org/10.1016/j.jnt.2020.06.015, preprint available as \href{https://arxiv.org/abs/1810.11941}{arXiv:1810.11941}.









\bibitem[AEK]{AEK} J. Arthur, D. Ellwood and R. Kottwitz \emph{Harmonic analysis, the trace formula, and Shimura varieties}, Proceedings of the Clay Mathematics Institute 2003, available as \href{http://www.claymath.org/library/proceedings/cmip04.pdf}{http://www.claymath.org/library/proceedings/cmip04.pdf}



\bibitem[BaPa]{BaPa} V. Balaji and A.J. Parameswaran, \emph{Semistable principal bundles-II (in positive
characteristics)}, Transformation Groups l8(1) (2003) 3--36, MR 1959761.

\bibitem[BaSe]{BaSe} V. Balaji, C.S. Seshadri, \emph{Semistable principal bundles. I (Characteristic zero)}, in: Special issue in celebration of
Claudio Procesi’s 60th birthday, J. Algebra 258 (2002) 321--347.




\bibitem[BK]{BK} S.~ Bloch, K.~Kato. \emph{p-adic \'etale cohomology}. Inst. Hautes \'Etudes Sci. Publ. Math., (63):107--152, 1986.






\bibitem[BoL\"u]{BL} S. Bosch, W. L\"utkebohmert, \emph{Formal and rigid geometry. II}. Flattening techniques, Math. Ann. 296 (1993), 403--429





\bibitem[BeLa]{B-L} A.\ Beauville and Y.\ Laszlo: \emph{Un lemme de descente}, Comptes Rendus Acad. Sci. Paris, vol. {\bfseries 320}, s\'erie I (1995), 335--340; also available as \href{http://math1.unice.fr/~beauvill/pubs/descente.pdf}{http:/\!/math1.unice.fr/$\sim$beauvill/pubs/bibli.html}.


\bibitem[Beh]{Beh} K.\ Behrend: \emph{The lefschetz trace formula for the moduli stack of principal bundles}, PhD thesis, University of California, Berkeley, 1991; available at \href{http://www.math.ubc.ca/~behrend/thesis.pdf}{http:/\!/www.math.ubc.ca/$\sim$behrend/}.

\bibitem[BeDh]{BeDh} K. Behrend and A. Dhillon: \emph{On the motivic class of the stack of bundles}. Adv. Math., 212(2):617--644, 2007.


\bibitem[BeDr]{B-D} A.~Beilinson, V.~Drinfeld: \emph{Quantization of Hitchin's integrable system and Hecke eigensheaves}, preprint on \href{http://www.math.uchicago.edu/~mitya/langlands/hitchin/BD-hitchin.pdf}{http:/\!/www.math.uchicago.edu/$\sim$mitya/langlands.html}.



\bibitem[BGM]{BGM} I. Biswas, T. L. Gomez, and V. Munoz, \emph{Automorphisms of moduli spaces of vector bundles over a curve}, Expo. Math. 31 (2013), no. 1, 73--86. available as \href{http://arxiv.org/abs/1202.2961}{arXiv:1202.2961}








\bibitem[BLR]{BLR} S.~Bosch, W.~L\"utkebohmert, M.~Raynaud: \emph{N\'eron models}, Ergebnisse der Mathematik und ihrer Grenzgebiete (3) {\bfseries 21}, Springer-Verlag, Berlin, 1990.


\bibitem[Br1]{Brion1} M. Brion, \emph{On automorphisms and endomorphisms of projective varieties}, in: Automorphisms
in birational and affine geometry, 59--81, Springer Proc. Math. Stat. 79, Springer, Cham, 2014.



\bibitem[Br2]{Brion2} M.~Brion \emph{On algebraic semigroups and monoids}. In Algebraic Monoids, Group Embeddings and Algebraic Combinatorics, vol. 71 of Fields Institute Communication Series, pp. 1-55. Springer, 2014. 

\Verkuerzung
{
}
{}



\Verkuerzung
{
}
{}


\bibitem[CD]{CD} D.-C. Cisinski and F. D\'eglise, \emph{Integral mixed motives in equal characteristic}, Doc. Math., Extra Volume: Alexander
S. Merkurjev’s Sixtieth Birthday (2015), pp. 145-194

\bibitem[Co]{Conrad} B. Conrad, Keel-mori theorem via stacks, available at http://www.math.stanford.edu/
~bdconrad/papers/coarsespace.pdf (2005).







\bibitem[Dri]{Drinfeld1} V.\ G.\ Drinfeld, \emph{Moduli varieties of $F$-sheaves}, Func.\ Anal.\ and Appl.\ {\bfseries 21} (1987), 107--122.




\bibitem[Fal]{Faltings03} G.~Faltings: \emph{Algebraic loop groups and moduli spaces  of bundles}, J.\ Eur.\ Math.\ Soc.\ {\bfseries 5} (2003), no.~1, 41--68. 


\bibitem[EGA]{EGA} A.~Grothendieck: \emph{{\'E}lements de G{\'e}o\-m{\'e}trie Alg{\'e}\-brique}, Publ.\ Math.\ IHES, 1960--1967; see also Grundlehren {\bfseries 166}, Springer-Verlag, Berlin etc.\ 1971.



\bibitem[Gro]{Gro} A.~Grothendieck : \emph{Techniques de construction et th\'eor\`emes d'existence en g\'eom\'etrie alg\'ebrique IV : les sch\'emas de Hilbert}, S\'eminaire Bourbaki 5 (1960--1961), Expos\'e No. 221.


\bibitem[HLP]{HLP} D. Halpern-Leistner and A. Preygel. \emph{Mapping stacks and categorical notions of properness}; 2014, arXiv:1402.3204v1 


\bibitem[Kol]{Kol} J.~Koll\'ar : \emph{Rational curves on algebraic varieties}, Ergeb. Math. Grenzgeb. (3) 32, Springer-Verlag, Berlin, 1996 










\bibitem[H-K]{H-K} 
A.\ Huber , B.\ Kahn: \emph{The slice filtration and mixed Tate motives}, \emph{Compos.\ Math.}, 142(4): 907-936 (2006).








\Verkuerzung
{
}
{}

\bibitem[Hei]{Heinloth} J.\ Heinloth: \emph{Uniformization of $\CG$-bundles}, Math.\ Ann.\ {\bfseries 347} (2010), 499--528; also available as \href{http://arxiv.org/abs/0711.4450}{arXiv:0711.4450}.


\bibitem[HeiSch]{HeiSch} J.\ Heinloth, A.\ Schmitt \emph{The cohomology ring of moduli stacks of principal bundles over curves.} Documenta Mathematica, Vol. 15, p. 423--488, 2010.


\bibitem[HoLe]{Vicky-Pepin} V. Hoskins and S. Pepin Lehalleur: \emph{On the Voevodsky motive of the moduli stack of vector bundles over a
curve}. arxiv: 1711.11072, 2017.
 




\bibitem[Kel]{Kel} S. Kelly, \emph{Triangulated categories of motives in positive characteristic}, PhD thesis of Université
Paris 13 and of the Australian National University, arXiv:1305.5349v2, 2012.




\bibitem[Kre]{Kre}
A. Kresch, \emph{Cycle groups for Artin stacks}, Invent. Math. 138 (1999), 495--536. 



\Verkuerzung
{
}
{}


\bibitem[Laf]{VLaff} V.~Lafforgue. \emph{Chtoucas pour les groupes r\'eductifs et param\'etrisation de Langlands globale.} Preprint \href{http://arxiv.org/abs/1209.5352}{arXiv:1209.5352}(2012).








\bibitem[MVW]{MVW} C. Mazza, V. Voevodsky, C. A. Weibel, \emph{Lecture notes on motivic cohomology},
Clay mathematics monographs, v.2., (2006).


\bibitem[Lau]{Lau} E.~ Lau, \emph{On generalised D-shtukas}, Dissertation, Rheinische Friedrich-Wilhelms Universit\"at Bonn, 2004. Bonner Mathematische Schriften 369. Universit\"at Bonn, Mathematisches Institut, Bonn (2004), available at the address http://www.math.uni-bielefeld.de/~lau/publ.html






\bibitem[Mil]{Milne} J. S. Milne, \emph{\'Etale cohomology}, Princeton Math. Series 33, Princeton Univ. Press, Princeton, N.J., 1980.




\bibitem[Ol1]{Ols1} M.~Olsson \emph{On proper coverings of  Artin stacks}, Adv.~ Math. 198 (2005), no. 1, 93--106.


\bibitem[Ol2]{Ols2}  M.~Olsson \emph{Hom--stacks and restriction of scalars}. Duke Math. J. 134 (2006), 139--164.

\bibitem[PoYu]{PY}
 M.~Porta and T.~ Y.~ Yu, \emph{Higher analytic stacks and GAGA theorems}, preprint 2014, http://arxiv.org/abs/1412.5166.



\bibitem[PaRa]{PR2} G.~Pappas, M.~Rapoport: \emph{Twisted loop groups and their affine flag varieties}, Advances in Math.\ {\bfseries 219} (2008), 118--198; also available as \href{http://arxiv.org/abs/math/0607130}{arXiv:math/0607130}.








\bibitem[HaRy]{HR} J.~ Hall, D. Rydh, \emph{Coherent Tannaka duality and algebraicity of Hom-stacks}, Algebra Number Theory 13 (2019), no. 7, 1633--1675.






\bibitem[Ric]{Richarz} T.~Richarz: \emph{Schubert varieties in twisted affine flag varieties and local models}, J.~Algebra {\bfseries 375} (2013), 121--147; also available as \href{http://arxiv.org/abs/1011.5416}{arXiv:1011.5416}.






\bibitem[Stacks]{Stacks} A. John de Jong (maintainer), The Stacks project: https://stacks.math.columbia.edu







\Verkuerzung
{
}
{}

\bibitem[Var]{Var} Y.\ Varshavsky: \emph{Moduli spaces of principal $F$-bundles}, Selecta Math.\ (N.S.) {\bfseries 10} (2004),  no.\ 1, 131--166; also available as \href{http://arxiv.org/abs/math/0205130}{arXiv:math/0205130}.

\bibitem[Vis]{Vis} A. Vistoli. \emph{Intersection theory on algebraic stacks and on their moduli spaces}, Inventiones mathematicae, 97; 613--670, 1989.


\bibitem[VSF]{VSF}
V.\ Voevodsly, A.\ Suslin, E.M.\ Friedlander: \emph{Cyles, Transfers, and Motivic Homology Theories}, Princeton university press (2000).



\bibitem[Wan]{Wang} J.~Wang: \emph{The moduli stack of $G$-bundles}, preprint on \href{http://arxiv.org/abs/1104.4828}{arXiv:1104.4828}.



\end{thebibliography}
}
\end{document}